\documentclass[a4paper, 11pt]{article}

\usepackage{amsmath, amsthm, amsfonts, hyperref}

\usepackage{enumitem}
\numberwithin{equation}{section}

\renewcommand{\Re}{\operatorname{Re}}
\renewcommand{\Im}{\operatorname{Im}}

\newcommand{\realR}{\mathbb{R}}

\DeclareMathOperator{\Tr}{Tr}
\DeclareMathOperator{\diag}{diag}

\newcommand{\gfrac}[2]{\genfrac{}{}{0pt}{}{#1}{#2}}
\newcommand{\MeijerG}[5]{G^{#1}_{#2} \left( \left. \gfrac{#3}{#4} \right\rvert {#5} \right)}
\newcommand{\hyperF}[5]{\,_{#1}F_{#2} \left( \left. \gfrac{#3}{#4} \right\rvert {#5} \right)}
\DeclareMathOperator{\kernel}{ker}
\DeclareMathOperator{\rank}{rank}
\newtheorem{theorem}{Theorem}[section]
\newtheorem{corollary}[theorem]{Corollary}
\newtheorem{cor}[theorem]{Corollary}

\newtheorem{lemma}[theorem]{Lemma}
\newtheorem{proposition}[theorem]{Proposition}
\theoremstyle{remark}
\newtheorem{rmk}[theorem]{Remark}

\newtheorem{assumption}[theorem]{Assumption}

\newcommand{\Andreief}{Andr\'{e}ief}
\newcommand{\Pineiro}{Pi\~{n}eiro}

\title{Correlation kernels for sums and products of random matrices}
\author{Tom Claeys\thanks{ Universit\'{e} catholique de Louvain, Chemin du cyclotron 2, B-1348 Louvain-La-Neuve, Belgium, \href{mailto:tom.claeys@uclouvain.be}{\nolinkurl{tom.claeys@uclouvain.be}} \newline
  Supported by the European Research Council under the European Union's Seventh Framework Programme (FP/2007/2013)/ ERC Grant Agreement n.~307074 and by the Belgian Interuniversity Attraction Pole P07/18} \and Arno B.~J.~Kuijlaars\thanks{KU Leuven, Department of Mathematics, Celestijnenlaan 200B Box 2400, 3001 Leuven, Belgium, \href{mailto:arno.kuijlaars@wis.kuleuven.be}{\nolinkurl{arno.kuijlaars@wis.kuleuven.be}} \newline
  Supported by KU Leuven Research Grant OT/12/073, the Belgian Interuniversity Attraction Pole P07/18, and FWO Flanders projects G.0641.11 and G.0934.13} \and Dong Wang\thanks{National University of Singapore, Department of Mathematics, Singapore, 119076, \href{mailto:matwd@nus.edu.sg}{\nolinkurl{matwd@nus.edu.sg}} \newline
    Supported partially by the startup grant R-146-000-164-133}}

\begin{document}

\maketitle
\begin{abstract}
  Let $X$ be a random matrix whose squared singular value density is a polynomial ensemble. We derive double contour integral formulas for the correlation kernels of the squared singular values of $GX$ and $TX$, where $G$ is a complex Ginibre matrix and $T$ is a truncated unitary matrix. We also consider the product of $X$ and several complex Ginibre/truncated unitary matrices. As an application, we derive the precise condition for the squared singular values of the product of several truncated unitary matrices to follow a polynomial ensemble. We also consider the sum $H + M$ where $H$ is a GUE matrix and $M$ is a random matrix whose eigenvalue density is a polynomial ensemble. We show that the eigenvalues of $H + M$ follow a polynomial ensemble whose correlation kernel can be expressed as a double contour integral. As an application, we point out a connection to the two-matrix model. 
\end{abstract}

\section{Introduction}

Eigenvalues and singular values of sums and products of random matrices have been studied extensively over the last decade. The focus of the research has been the global density of the eigenvalues and singular values of the sum and product of large random matrices. As the dimension of random matrices tends to infinity, the macroscopic eigenvalue density of sums and products of random matrices can be described under fairly general conditions using free probability techniques \cite{Nica-Speicher06}. However, the free probability method can only give global asymptotic results. For exact correlation functions of the eigenvalues and singular values of the sums and products of random matrices, until recently there were only a limited number of cases in which results were known. Notable cases are the Gaussian unitary ensemble with external source \cite{Brezin-Hikami96}, which is the sum of a Gaussian unitary matrix (GUE) with a deterministic matrix; and the complex Wishart ensemble \cite{Baik-Ben_Arous-Peche05}, which is the product of a complex Ginibre matrix with a deterministic matrix. Another interesting result is known for the sum of complex null Wishart matrices, in the form of the multiple Laguerre minor process with a fixed time \cite{Adler-van_Moerbeke-Wang11}. It is remarkable that all the examples mentioned above are \emph{polynomial ensembles}, a special kind of determinantal point processes, and that their correlation kernels all have double contour integral representations.

Recently, the correlation functions of the (squared) singular values are found for more types of products of random matrices, for instance for products of complex Ginibre matrices and products of truncated unitary Haar distributed matrices, as well as their inverses, see \cite{Akemann-Ipsen-Kieburg13}, \cite{Kuijlaars-Zhang14}, \cite{Forrester14}, \cite{Kuijlaars-Stivigny14}, \cite{Kieburg-Kuijlaars-Stivigny15}, and \cite{Forrester-Liu15}. For a comprehensive survey of the current developments, see \cite{Akemann-Ipsen15} and references therein. In all the cases listed above, the squared singular value densities are polynomial ensembles, and the correlation kernels can be expressed as double contour integrals.

Polynomial ensembles arise naturally in the study of probability densities for eigenvalues and singular values of random matrices. Their structure is preserved under certain operations on random matrices, such as multiplication by a complex Ginibre matrix \cite{Kuijlaars-Stivigny14} and multiplication by a truncated unitary Haar distributed matrix \cite{Kieburg-Kuijlaars-Stivigny15}. Those random matrix operations thus induce transformations of polynomial ensembles \cite{Kuijlaars15}. In this paper, we obtain explicit transformation formulas for the correlation kernels for a number of such transformations.

One feature of our results is that the correlation kernels obtained in this paper are all in the double contour integral form. Double contour integral formulas are very suitable to derive asymptotic results on local statistics, by the classical saddle point method. Numerous applications of double contour formulas can be found in literature. Two recent examples are \cite{Kuijlaars-Zhang14} for the discovery of a family of new hard-edge universality and \cite{Liu-Wang-Zhang14} for the proof of bulk and soft-edge universality.

\medskip 

A polynomial ensemble \cite{Kuijlaars-Stivigny14} is a probability density function for $n$ particles $x_1, \ldots, x_n$ on the real
line of the form
\begin{equation} \label{f-pol-ensemble} 
	\frac{1}{Z_n} \Delta_n(x) \det \left[ f_{k-1}(x_j) \right]_{j,k=1}^n 
	\end{equation}
with given functions $f_0, \ldots, f_{n-1}$.
Here we use  
\begin{equation}
  \Delta_n(x) = \prod_{1\leq j<k\leq n} (x_k-x_j) = \det \left[ x_j^{k-1}\right]_{j,k=1}^n
\end{equation} 
to denote the Vandermonde determinant for the $n$-tuple $x = (x_1, \ldots, x_n) \in \mathbb R^n$, and $Z_n$ is a normalization constant.
It is a polynomial ensemble on $E \subset \mathbb R$ if 
$f_k(x) = 0$ for every $x \in \mathbb R \setminus E$
and every $k=0, \ldots, n-1$, since then the particles are in $E$
(with probability one).

The polynomial ensemble \eqref{f-pol-ensemble} is a determinantal point process and its correlation kernel can be written in the form
\begin{equation} \label{eq:Kn} 
  k_n(x,y) = \sum_{j=0}^{n-1} p_j(x) q_j(y), 
\end{equation}
where $p_j$ is a polynomial of degree $j$ and $q_j$ belongs to the linear
span of $f_0, \ldots, f_{n-1}$ for every $j=0, \ldots, n-1$, in such a way that they satisfy the biorthogonality conditions
\begin{equation} 
	\label{eq:pqbiort} \int_{-\infty}^{\infty} p_j(x) q_k(x) dx = \delta_{j,k},
	\qquad \text{for } j, k = 0, \ldots, n-1.
\end{equation}
We call such a set of $p_j$'s and $q_j$'s a {\em biorthogonal system} associated with 
the polynomial ensemble \eqref{f-pol-ensemble}. We sometimes refer to the $q_j$'s as {\em dual functions}.
The $p_j$'s and $q_j$'s are not unique, but the sum $\sum_{j=0}^{n-1} p_j(x) q_j(y)$ in \eqref{eq:Kn} is.
If we require that $p_j$ is a monic polynomial, and that $q_j$ is in the linear
span of $f_0, \ldots, f_j$ for every $j$, then the biorthogonal system is unique.
However, such a biorthogonal system does not always exist, and, if it exists, it may not 
be the most convenient one to work with. The unique monic polynomial $p_n$ of degree $n$ which is orthogonal to $f_0,\ldots, f_{n-1}$, i.e.\ $\int_{-\infty}^{\infty}p_n(x)f_k(x)dx=0$ for $k=0,\ldots, n-1$, is the average characteristic polynomial 
\begin{equation} \label{eq:p_i_as_char_func}
  p_n(x) = \mathbb E \left( \prod^n_{j = 1} (x - x_j) \right),
\end{equation} 
where $\mathbb E$ denotes the average over particles $x_1,\ldots, x_n$ in the polynomial ensemble \eqref{f-pol-ensemble}.

In \cite{Kuijlaars15} a number of transformations which preserve the structure of a polynomial ensemble 
are given with transformation formulas for the functions $f_k$ in \eqref{f-pol-ensemble}, but without explicit formulas for 
the correlation kernels. 
In this paper, we give transformation formulas for the correlation kernels, biorthogonal systems and 
average characteristic polynomials for three such transformations. The polynomial ensembles which we 
consider are joint probability densities for:
\begin{enumerate}[label=(\roman*)]
\item \label{enu:case_i}
  the eigenvalues of the sum $H+M$ of a GUE matrix $H$ with a Hermitian random matrix $M$,
\item \label{enu:case_ii}
  the squared singular values of the product $GX$ of a complex Ginibre matrix $G$ with a random matrix $X$,
\item \label{enu:case_iii}
  the squared singular values of the product $TX$ of a truncation of a Haar distributed
	unitary matrix with a random matrix $X$.
\end{enumerate}

Here and throughout the paper, we understand $M$ and $X$ as follows.

\begin{assumption} \label{defn:M_matrix}
      $M$ is an $n \times n$ Hermitian random matrix whose eigenvalues $x_1, \dotsc, x_n$ follow a polynomial ensemble 
			expressed by \eqref{f-pol-ensemble}. We let $p_k(x)$ and $q_k(y)$ ($k = 0, 1, \dotsc, n - 1$) be 
			a biorthogonal system, $k_n(x, y)$ be the correlation kernel, and $p_n(x)$ be the average characteristic polynomial 
			for the polynomial ensemble.
    \end{assumption}

\begin{assumption} \label{defn:X_matrix}
      $X$ is an $\ell \times n$ ($\ell \geq n$) complex matrix whose squared singular values 
			$x_1, \dotsc, x_n$, that is, eigenvalues of $X^* X$, follow a polynomial ensemble 
			expressed by \eqref{f-pol-ensemble}. We let $p_k(x)$ and $q_k(y)$ ($k = 0, 1, \dotsc, n - 1$) 
			be the biorthogonal system, $k_n(x, y)$ be the correlation kernel, and $p_n(x)$ 
			be the average characteristic polynomial for the polynomial ensemble.
\end{assumption}
Note that although notations $f_k$, $p_k$, $q_k$, $k_n$ and $p_n$ have different meanings in Assumptions
\ref{defn:M_matrix} and \ref{defn:X_matrix}, since $M$ and $X$ will not appear together in our paper, no confusion will occur.

In cases \ref{enu:case_ii} and \ref{enu:case_iii}, transformation formulas were obtained in \cite{Kuijlaars-Stivigny14} 
and  \cite{Kieburg-Kuijlaars-Stivigny15} for the joint probability densities of the squared singular values of 
$GX$ and $TX$ in terms of the density for $X$. We will prove a similar transformation formula in case \ref{enu:case_i}. 
Our main focus is on the correlation kernels. We will express, in case \ref{enu:case_i}, the eigenvalue 
correlation kernel of $H+M$ in terms of the eigenvalue correlation kernel of $M$, and in case \ref{enu:case_ii} 
(resp.\ case \ref{enu:case_iii}) the squared singular value correlation kernels of $GX$ (resp.\ $TX$), in terms of 
the squared singular value correlation kernel of $X$. In addition, we will give transformation formulas for biorthogonal 
systems and for the average characteristic polynomials.

\section{Statement of main results}

\subsection{Addition of a GUE matrix}

Up to an overall scale factor, a GUE matrix is a Hermitian random matrix $H$ with diagonal entries in standard real normal distribution, upper-diagonal entries in standard complex normal distribution, and all upper-triangular and diagonal entries independent.
The joint probability density function for the eigenvalues $x_1,\ldots, x_n$ of $H$ is given by
\begin{equation} \label{eq:distr_GUE}
\frac{1}{Z_n}\Delta_n(x)^2\prod_{j=1}^ne^{-\frac{x_j^2}{2}}, \qquad Z_n = (2\pi)^{\frac{n}{2}} \prod^n_{k = 1} k!,
\end{equation}
which is the polynomial ensemble \eqref{f-pol-ensemble} in the case where $f_k(x)=x^ke^{-x^2/2}$. If we add a GUE matrix to a random matrix with eigenvalues in a polynomial ensemble, we are led to a transformed polynomial ensemble for the eigenvalues of the sum.

\begin{theorem}\label{theorem GUE jpdf}
  Let $H$ be an $n \times n$ GUE matrix and $M$ be the random matrix from
	Assumption \ref{defn:M_matrix}, independent of $H$. The density of the eigenvalues $y_1,\ldots, y_n$ of $H+M$ is given by
\begin{equation} \label{f-pol-ensemble 2} 
	\frac{1}{Z_n'} \Delta_n(y) \det \left[ F_{k-1}(y_j) \right]_{j,k=1}^n,  
	\end{equation}
	for some constant $Z_n'$,
	where
	\begin{equation}\label{def F GUE}
    F_k(y) = \int_{-\infty}^{\infty} f_k(t) e^{-\frac{1}{2}(y-t)^2} dt,\qquad k=0,\ldots, n-1.
  \end{equation}
  \end{theorem}

We will prove Theorem \ref{theorem GUE jpdf} in Section \ref{section: proof GUE jpdf}.

\begin{rmk}
The functions $F_k$ in the transformed polynomial ensemble \eqref{f-pol-ensemble 2} are convolutions of the functions $f_k$ from the original polynomial ensemble \eqref{f-pol-ensemble} with the Gaussian density $e^{-x^2/2}$.  
\end{rmk} 
 
We now describe the effect of this transformation on the biorthogonal system, on the correlation kernel, and on the average characteristic polynomial.
Since the eigenvalue density of $H+M$ is the polynomial ensemble \eqref{f-pol-ensemble 2}, the associated correlation kernel $K_n$ can be written in the form
\begin{equation} \label{eq:Kn hat} 
	K_n(x,y) = \sum_{k=0}^{n-1} P_k(x) Q_k(y), 
	\end{equation}
where $P_k$ is a polynomial of degree $k$ and $Q_k$ is in the linear
span of $F_0, \ldots, F_{n-1}$,  such that
\begin{equation} 
	\label{eq:PQbiort} \int_{-\infty}^{\infty} P_j(x) Q_k(x) dx = \delta_{j,k},
	\qquad \text{for } j,k = 0, \ldots, n-1.
\end{equation}
We express the quantities $K_n$, $P_k$, and $Q_k$ in the transformed polynomial ensemble in terms of their counterparts $k_n$, $p_k$, and $q_k$ in the original polynomial ensemble for $M$.

\begin{theorem}\label{theorem: kernel GUE}
 Let $H$ be an $n \times n$ GUE random matrix  and $M$ be the random matrix 
from Assumption  \ref{defn:M_matrix}, independent of $H$, as in Theorem \ref{theorem GUE jpdf}. Then, 
 \begin{enumerate}[label=(\alph*)]
 \item \label{enu:theorem: kernel GUE_a}
   a biorthogonal system in the transformed polynomial ensemble defined by \eqref{f-pol-ensemble 2} and \eqref{def F GUE} is given by
  \begin{align}
    P_k(x) = {}& \frac{1}{\sqrt{2\pi} i} \int_{-i\infty}^{+i\infty} p_k(s) e^{\frac{1}{2}(x-s)^2} ds, \label{Pk} \\
    Q_k(y) = {}& \frac{1}{\sqrt{2\pi}} \int_{-\infty}^{\infty} q_k(t) e^{-\frac{1}{2}(y-t)^2}dt, \label{Qk}
  \end{align}
  \item \label{enu:theorem: kernel GUE_b}
    the correlation kernel for the eigenvalues of $H+M$ is given by
    \begin{equation}\label{integral transform}
      K_n(x, y)=\frac{1}{2\pi i} \int_{-i\infty}^{+i\infty}ds \int_{-\infty}^{\infty}dt\, k_n \left(s,t\right)e^{\frac{1}{2}((x- s)^2-(y- t)^2)},
    \end{equation}
  \item \label{enu:theorem: kernel GUE_c}
    the average characteristic polynomial $P_n$ of $H+M$ is given by \eqref{Pk} with $k=n$.
  \end{enumerate}
\end{theorem}

The proof of Theorem \ref{theorem: kernel GUE} will be given in Section \ref{section: proof GUE}.

\begin{rmk}
In the integral formula \eqref{integral transform}, we use an extension of the kernel $k_n(s,t)$ for $s$ in the complex plane. We can do this because $k_n(s,t)$ is polynomial in $s$ for any $t\in\mathbb R$.
\end{rmk}
\begin{rmk}\label{rmk: Weierstrass}
The integral formulas \eqref{Qk} and \eqref{def F GUE} for $Q_k$ and $F_k$ can be recognized as Weierstrass transforms \cite[Chapter VIII]{Hirschman-Widder55} of $q_k$ and $f_k$ (up to the prefactors; the usual definition of the Weierstrass transform also has $1/4$ instead of $1/2$ in the exponent); $P_k$ given in \eqref{Pk} is  the inverse Weierstrass transform of $p_k$, up to the prefactor, which is chosen such that $P_k$ is monic if $p_k$ is monic. 
\end{rmk}

As an application of Theorem \ref{theorem: kernel GUE}, we derive in Section \ref{subsubsec:GUE+1MM} a double contour integral formula \eqref{kernelGUE+UECauchy} for the correlation kernel of the eigenvalues of the sum of a GUE matrix and an arbitrary random unitary invariant matrix. This formula is suitable for the analysis of limiting local statistics, and we investigate the analytic consequences in a separate paper.

\subsection{Multiplication with a Ginibre matrix}

A complex Ginibre matrix is a rectangular matrix whose entries are independent and identically distributed in the standard complex normal distribution. We define the complex Ginibre matrix $G$ as follows.
\begin{assumption} \label{defn:G_matrix}
      Let $\nu \geq 0$ and $\ell \geq n$. $G$ is an $(n + \nu) \times \ell$ complex Ginibre matrix.
\end{assumption}
In the case $\ell = n$, the joint probability density for the squared singular values $x_1,\ldots, x_n$ of $G$ is given by
\begin{equation}
  \frac{1}{Z_n}\Delta_n(x)^2\prod_{j=1}^n x_j^\nu e^{-x_j},\qquad x_j\geq 0, \qquad Z_n = \prod^n_{k = 1} k! \Gamma(k + \nu).
\end{equation}
This is a polynomial ensemble on $[0,\infty)$ defined by the functions $f_k(x)=x^{\nu+k} e^{-x}$, $k = 0, 1, \dotsc, n - 1$.

Left multiplication of a random matrix by a complex Ginibre matrix induces a transformation of polynomial ensembles on $[0,\infty)$: the following result was proved by Kuijlaars and Stivigny \cite{Kuijlaars-Stivigny14} (based on ideas from \cite{Akemann-Ipsen-Kieburg13, Akemann-Kieburg-Wei13}).
\begin{proposition} \label{GinibreTransformation}
  Let $G$ be the complex Ginibre matrix from Assumption \ref{defn:G_matrix} and $X$ be the random matrix 
	from Assumption \ref{defn:X_matrix}, independent of $G$. Then the squared singular values $y_1, \ldots, y_n$ of $Y = GX$ follow a polynomial ensemble on $[0,\infty)$ given by
\begin{equation} \label{F-pol-ensemble} 
	\frac{1}{{Z}_n'} \Delta_n(y) \, \det \left[ F_{k-1}(y_j) \right]_{j,k=1}^n,  
	\end{equation}
	for some constant $Z_n'$,
where 
\begin{equation} \label{Fk-definition} 
	F_k(y) = \int_0^{\infty} t^{\nu} e^{-t} f_k \left( \frac{y}{t} \right) \frac{dt}{t}, \qquad y > 0.
	\end{equation}
\end{proposition}

We complement the above result with transformation formulas for the biorthogonal system, for the correlation kernel, and for the average characteristic polynomial.

\begin{theorem} \label{GinibreTransformationKernel}
  Let $G$ be the complex Ginibre matrix from Assumption \ref{defn:G_matrix} and $X$ be the random matrix 
	from Assumption \ref{defn:X_matrix}, independent of $G$, as in Proposition \ref{GinibreTransformation}. Then,
  \begin{enumerate}[label=(\alph*)]
  \item \label{enu:GinibreTransformationKernel_a}
    a biorthogonal system in the transformed polynomial ensemble defined by \eqref{F-pol-ensemble} and \eqref{Fk-definition} is given by
    \begin{align}
      P_k(x) = {}& \frac{1}{2\pi i} \oint_{\Sigma} s^{-\nu} e^{s} p_k\left(\frac{x}{s}\right) \frac{ds}{s}, \label{eq:biorthogonal_pd_P} \\
      Q_k(y) = {}& \int_{0}^{\infty} t^{\nu} e^{-t} q_k\left(\frac{y}{t}\right) \frac{dt}{t}, \label{eq:biorthogonal_pd_Q}
    \end{align}
    where $\Sigma$ is a simple counter-clockwise oriented contour encircling the origin,
  \item \label{enu:GinibreTransformationKernel_b}
    the correlation kernel $K_n$ for the squared singular values of $GX$ is given by
    \begin{equation}\label{integral transform product}
      K_n(x, y)  = \frac{1}{2\pi i} \oint_\Sigma\frac{ds}{s} \int_{0}^{\infty}\frac{dt}{t} \,  k_n\left(\frac{x}{s},\frac{y}{t}\right)
			\left(\frac{t}{s}\right)^\nu e^{s-t},
    \end{equation}
  \item \label{enu:GinibreTransformationKernel_c}
    the average characteristic polynomial $P_n$ of $X^*G^*GX$ is equal to
		\begin{equation} \label{eq:biorthogonal_pd_Pn}
		P_n(x) =  \frac{(n+\nu)!}{2\pi i} \oint_{\Sigma} s^{-\nu} e^{s} p_n\left(\frac{x}{s}\right) \frac{ds}{s}. \end{equation}
  \end{enumerate}
\end{theorem}	  

The proof of Theorem \ref{GinibreTransformationKernel} will be given in Section \ref{section: proof Ginibre}.

\begin{rmk}\label{rmk: Mellin Gamma}
The functions $Q_k$ and $F_k$ can be recognized as Mellin convolutions of $f_k$ and $q_k$ with the 
Gamma density $t^{\nu} e^{-t}$. We will show that $P_k$ defined in \eqref{eq:biorthogonal_pd_P} can also be written in an alternative way: 
if $p_k(x)=\sum_{j=0}^ka_jx^j$, then we have
\begin{equation}\label{Pk Hadamard}
  P_k(x) = \sum_{j=0}^k \frac{a_j}{(j+\nu)!}x^j.
\end{equation}
Thus the factor $(n + \nu)!$ in \eqref{eq:biorthogonal_pd_Pn} makes $P_n(x)$ monic.. The transformation formula \eqref{Pk Hadamard} for $P_k$ is reminiscent of similar transforms in \cite[Formulas (3.16) and (3.43)]{Forrester-Liu14}. We can also write \eqref{integral transform product} as
\begin{equation}\label{integral transform product alternative}
  K_n(x, y)  = \frac{1}{2\pi i}\left(\frac{y}{x}\right)^\nu \oint_\Sigma\frac{ds}{s} 
	\int_{0}^{\infty}\frac{dt}{t} \, k_n(s,t)\left(\frac{s}{t}\right)^\nu e^{\frac{x}{s}-\frac{y}{t}}.
\end{equation}
\end{rmk}

\subsection{Multiplication with a truncated unitary matrix}

The truncated unitary matrix $T$ is defined as follows.
\begin{assumption} \label{defn:T_matrix}
        Let $\nu \geq 1$, $m \geq \ell \geq n$ and $\mu = m - n - \nu \geq 1$. Let $U$ be an 
				$m \times m$ Haar distributed random unitary matrix. $T$ is the $(n + \nu) \times \ell$ truncation of $U$.
\end{assumption}
In the case where $\ell = n$ and $\mu \geq n$, the joint probability density for the squared 
singular values $x_1,\ldots, x_n$ of $T$ is given by \cite[Section 3.8.3]{Forrester10}
\begin{multline} \label{Jacobi-ensemble}
  \frac{1}{Z_n} \Delta_n(x)^2 \prod^n_{j = 1} x^{\nu}_j (1 - x_j)^{m - 2n - \nu}, \qquad 0 \leq x_j \leq 1, \\
  Z_n = \prod^n_{k = 1} \frac{k! \Gamma(k + \nu) \Gamma( k + m - 2n - \nu)}{\Gamma(k + m - n)}.
\end{multline}

The following result is due to Kieburg, Kuijlaars and Stivigny \cite{Kieburg-Kuijlaars-Stivigny15}.

\begin{proposition} \label{TruncationTransformation}
  Let $T$ be the truncated unitary matrix from Assumption \ref{defn:T_matrix} and $X$ be the random matrix 
	from Assumption \ref{defn:X_matrix}, independent of $T$. Then the squared singular values $y_1, \ldots, y_n$ of $Y = TX$ follow a polynomial ensemble on $[0,\infty)$ given by
  \begin{equation} \label{F-pol-ensemble2} 
    \frac{1}{{Z}_n'} \, \Delta_n(y) \, \det \left[ F_{k-1}(y_j) \right]_{j,k=1}^n, 
  \end{equation}
  for some constant $Z_n'$,
  where 
  \begin{equation} \label{Fk-definition2} 
    F_k(y) = \int_0^1 t^{\nu} (1-t)^{\mu-1} f_k \left( \frac{y}{t} \right) \frac{dt}{t}, \qquad y > 0.
  \end{equation}
\end{proposition}
We have the counterpart of Theorem \ref{GinibreTransformationKernel} as follows.
\begin{theorem} \label{TruncationTransformationKernel}
  Let $T$ be the truncated unitary matrix from Assumption \ref{defn:T_matrix} and $X$ be the random matrix 
  from Assumption \ref{defn:X_matrix}, independent of $T$, as in Proposition \ref{TruncationTransformation}. Then, 
  \begin{enumerate}[label=(\alph*)]
  \item \label{enu:TruncationTransformationKernel_a}
    a biorthogonal system in the transformed polynomial ensemble defined by \eqref{F-pol-ensemble2} and \eqref{Fk-definition2} is given by
    \begin{align}
      P_k(x) & = \frac{\mu}{2\pi i} \oint_{\Sigma} s^{-\nu} (1-s)^{-\mu-1} p_k \left( \frac{x}{s} \right) \frac{ds}{s},  \label{def:Pk3} \\
      Q_k(y) & = \int_0^{1} t^{\nu} (1-t)^{\mu-1} q_k \left( \frac{y}{t} \right) \frac{dt}{t}, \label{def:Qk3}
    \end{align}
    where $\Sigma$ is a simple counter-clockwise contour around $0$ but not containing $1$,
  \item \label{enu:TruncationTransformationKernel_b}
    the correlation kernel for the squared singular values of $TX$ is given by
    \begin{equation} \label{integral transform trunc}
      {K}_n(x,y)  = 
      \frac{\mu}{2\pi i} \oint_{\Sigma} \frac{ds}{s} \int_0^1 \frac{dt}{t} 
      \left(\frac{t}{s} \right)^{\nu}  
      \, 
      k_n \left( \frac{x}{s}, \frac{y}{t} \right) \left(1- s \right)^{-\mu-1}  \left(1- t \right)^{\mu-1},
    \end{equation}
  \item \label{enu:TruncationTransformationKernel_c}
    the average characteristic polynomial $P_n$ of $X^*T^*TX$ is equal to
		\begin{align} \label{def:Pn3}
      P_n(x) & = \frac{\mu! (n+\nu)!}{2\pi i (n+\nu+\mu)!} \oint_{\Sigma} s^{-\nu} (1-s)^{-\mu-1} 
			p_n \left( \frac{x}{s} \right) \frac{ds}{s}.
					\end{align}
  \end{enumerate}
\end{theorem}
We will prove Theorem \ref{TruncationTransformationKernel} in Section \ref{section: proof truncated}.
\begin{rmk}\label{rmk: Mellin Beta}
  The functions $Q_k$ and $F_k$ are again Mellin convolutions of $q_k$ and $f_k$, but now with the Beta density $t^{\nu}(1-t)^{\mu-1}$ supported on $[0,1]$. We will show that $P_k$ in \eqref{def:Pk3} can alternatively be written as 
  \begin{equation}\label{Pk Hadamard trunc}
    P_k(x) = \frac{1}{(\mu-1)!} \sum_{j=0}^k \frac{(j+\nu+\mu)!}{(j+\nu)!}a_jx^j,
  \end{equation}
  if $p_k(x)=\sum_{j=0}^ka_jx^j$. Thus the factor $\mu! (n+\nu)!/(n+\nu+\mu)!$ in \eqref{def:Pn3} makes $P_n(x)$ monic. The formula \eqref{Pk Hadamard trunc} is analogous to \eqref{Pk Hadamard}, and is reminiscent of similar formulas in \cite{Forrester-Liu14}, see Remark \ref{rmk: Mellin Gamma}.
	We can also write \eqref{integral transform trunc} as
  \begin{equation}
    K_n(x, y) = \frac{\mu}{2\pi i} \left( \frac{y}{x} \right)^{\nu} \oint_{\Sigma'} \frac{ds}{s - x} \int^{\infty}_y \frac{dt}{t - y} k_n(s, t) \left( \frac{s}{t} \right)^{\nu + \mu} \left( \frac{t - y}{s - x} \right)^{\mu},
  \end{equation}
  where $\Sigma'$ is a simple counter-clockwise contour encircling both $0$ and $x$.
\end{rmk}

\subsection{Generalization}

Proposition \ref{GinibreTransformation}, Theorem \ref{GinibreTransformationKernel}, 
Proposition \ref{TruncationTransformation}, and Theorem \ref{TruncationTransformationKernel} have an obvious 
structural similarity. 
The new functions $F_k$ in \eqref{Fk-definition} and \eqref{Fk-definition2}
are given as Mellin convolutions of $f_k$ with a fixed function $\varphi$
where $\varphi(t) = t^{\nu} e^{-t}$ in the case of multiplication with
a complex Ginibre matrix as in Proposition \ref{GinibreTransformation} and
$\varphi(t) = t^{\nu} (1-t)^{\mu-1} \chi_{[0,1]}(t)$ in the case of multiplication
with a truncated unitary matrix as in Proposition \ref{TruncationTransformation}. 
The two Theorems \ref{GinibreTransformationKernel} and \ref{TruncationTransformationKernel} 
are obtained as  special cases of the following lemma which deals with
a general function $\varphi$. 

\begin{lemma} \label{Transformation-phi}
Let $p_0, \ldots, p_{n-1}$, $q_0, \ldots, q_{n-1}$ be a biorthogonal system for
the polynomial ensemble \eqref{f-pol-ensemble} on $[0,\infty)$. 
Suppose
\begin{equation} \label{Fphi-pol-ensemble} 
	\frac{1}{{Z}_n'} \Delta_n(y) \, \det \left[ F_{k-1}(y_j) \right]_{j,k=1}^n 
	\end{equation}
is a polynomial ensemble on $[0,\infty)$ with 
\begin{equation} \label{Fkphi-definition} 
	F_k(y) = \int_0^{\infty} \varphi(t) f_k \left( \frac{y}{t} \right) \frac{dt}{t}, \qquad y > 0,
	\end{equation}
where $\varphi : [0,\infty) \to [0,\infty)$ is a given non-negative
function with finite non-zero moments.
Define
\begin{equation} \label{eq:bj} 
	b_j = \left[ \int_0^{\infty} t^j \varphi(t) \, dt \right]^{-1} > 0,
	\qquad j=0,1, \ldots, n-1. \end{equation}
\begin{enumerate}[label=(\alph*)]
\item \label{enu:ransformation-phi_a}
  Then a biorthogonal system for \eqref{Fphi-pol-ensemble} is given 
  by polynomials $P_0, \ldots, \linebreak[4] P_{n-1}$,
  and functions $Q_0, \ldots, Q_{n-1}$ where for $k=0,1, \ldots, n-1$,
  \begin{equation} \label{eq:Pkphiseries} 
    P_k(x) = \sum_{j=0}^k  a_{j,k} b_j x^j \qquad \text{if} \quad p_k(x) = \sum_{j=0}^k a_{j,k} x^j 
  \end{equation}
  and
  \begin{equation} \label{eq:Qkphi} 
    Q_k(y) = \int_0^{\infty} \varphi(t) q_k \left( \frac{y}{t} \right) \frac{dt}{t}, \qquad y > 0.
  \end{equation}
\item \label{enu:ransformation-phi_b}
  Let $\psi$ be given by
  \begin{equation} \label{eq:psi}
    \psi(x) = \sum_{j=-\infty}^{\infty} b_j x^j 
  \end{equation}
  where $b_j$ is given by \eqref{eq:bj} for $j=0, \ldots, n - 1$ and otherwise arbitrary.
  Assume that the Laurent series \eqref{eq:psi} converges in the annulus $\{ r_{\psi} < |x| < R_{\psi} \}$.
  Then the polynomial $P_k$ has the alternative representation
  \begin{equation} \label{eq:Pkphiintegral} 
    P_k(x) = \frac{1}{2\pi i} \oint_{\Sigma} \psi(s) p_k \left(\frac{x}{s}\right) \frac{ds}{s}, 
  \end{equation}
  where $\Sigma$ is a closed, positive oriented contour lying in the annulus and encircling the origin once.
\item \label{enu:ransformation-phi_c}
 Under the same assumptions as in (b), if $k_n(x,y)$ is the correlation kernel for the polynomial
  ensemble \eqref{f-pol-ensemble}, then the correlation kernel $K_n(x,y)$
  for the transformed ensemble \eqref{Fphi-pol-ensemble} is given by
  \begin{equation} \label{eq:Knhatphi}
    K_n(x,y) = 
    \frac{1}{2\pi i} \oint_{\Sigma} \frac{ds}{s} \int_0^{\infty} \frac{dt}{t} 
    \psi(s) \varphi(t) k_n \left( \frac{x}{s}, \frac{y}{t} \right).
  \end{equation}
\end{enumerate}
\end{lemma}
We prove Lemma \ref{Transformation-phi} in Section \ref{sec:general_result}
and we obtain Theorems \ref{GinibreTransformationKernel} and \ref{TruncationTransformationKernel}
from it as easy consequences.

\subsection{Extensions}

Lemma \ref{Transformation-phi} can be applied to other situations as well, in particular
to products with a number of Ginibre matrices or a number of truncated unitary matrices. 
We first consider the product of a random matrix $X$ satisfying Assumption \ref{defn:X_matrix}
and several complex Ginibre matrices as follows.
\begin{assumption} \label{defn:G_matrices}
      Let $\nu_0,\nu_1, \dotsc, \nu_r$ be non-negative integers. Let
			$G_1, G_2, \dotsc, \linebreak[4] G_r$ be independent Ginibre matrices where $G_j$ is of size $(n + \nu_j) \times (n + \nu_{j - 1})$ for $j = 1,2, \dotsc, r$.
\end{assumption}
By applying 
Proposition \ref{GinibreTransformation} repeatedly, one easily sees that the squared singular values of $G_r G_{r - 1} \dotsm G_1 X$ follow a polynomial ensemble, as shown in \cite{Kuijlaars-Stivigny14}. 
It fits in the framework of Lemma \ref{Transformation-phi} with the function
\begin{equation} \label{phi-Mconvolution}
  \varphi = \varphi_r \ast \cdots \ast \varphi_1, 
\end{equation}
where $ \varphi_j(t) = t^{\nu_j} e^{-t}$ for $j=1, \ldots, r$, and
$\ast$ denotes the Mellin convolution.
Then $\varphi$ is a Meijer G-function,  $\varphi(t) = \MeijerG{r,0}{0,r}{-}{\nu_1, \ldots, \nu_r}{t}$,
and we use Meijer G-functions freely in our results concerning products with several matrices.  
For their definition and properties, see \cite{Luke69}, \cite{Boisvert-Clark-Lozier-Olver10},
and \cite{Beals-Szmigielski13}.

Then we have, from Lemma \ref{Transformation-phi}:
\begin{corollary} \label{cor:Ginibre_iterated}
  Let $G_1, \dotsc, G_r$ be the complex Ginibre matrices from Assumption \ref{defn:G_matrices} and let $X$ be the random matrix 
	from Assumption \ref{defn:X_matrix} with $\ell=n+\nu_0$, independent of $G_1, \dotsc, G_r$. Then the squared singular value density of the product $G_r \cdots G_1 X$ is a polynomial ensemble on $[0, \infty)$ with biorthogonal system consisting of polynomials
  \begin{equation}
    P_k(x) = \frac{1}{2\pi i} \oint_{\Sigma} \MeijerG{1,1}{1,r+1}{0}{0,-\nu_1, \ldots, -\nu_r}{-s} p_k \left( \frac{x}{s} \right) \frac{ds}{s}
  \end{equation}
  where $\Sigma$ is a simple closed counter-clockwise oriented contour around $0$,
  and dual functions
  \begin{equation}
    Q_k(y) = \int_0^{\infty} \MeijerG{r,0}{0,r}{-}{\nu_1, \ldots, \nu_r}{t} 
    q_k \left(\frac{y}{t}\right) \frac{dt}{t}.
  \end{equation}
  Furthermore, the correlation kernel for the transformed polynomial ensemble is
  \begin{multline}  \label{Knhat-rGinibre}
    K_n(x,y) = \frac{1}{2\pi i} \oint_{\Sigma} \frac{ds}{s} \int_0^{\infty} \frac{dt}{t}
    \MeijerG{1,1}{1,r+1}{0}{0, -\nu_1, \ldots, -\nu_r}{-s} \\
    \times \MeijerG{r,0}{0,r}{-}{\nu_1, \ldots, \nu_r}{t}
    k_n \left(\frac{x}{s}, \frac{y}{t} \right).
  \end{multline}
\end{corollary}

The proof of Corollary \ref{cor:Ginibre_iterated} will be given in Section \ref{section: corollary Ginibre}.

Similarly, we consider the product of the random matrix $X$ and several truncated unitary matrices as follows.
\begin{assumption} \label{defn:T_matrices}
     Let $\nu_0,\nu_1, \dotsc, \nu_r$ be non-negative integers. Let  
		$m_j \geq n + \nu_{j - 1}$ and $m_j - n - \nu_j = \mu_j \geq 1$ for all 
		$j = 1, 2, \dotsc, r$. Let $U_1, \dotsc, U_r$ be independent Haar distributed random unitary matrices 
		of size $m_j$ for $j=1, \ldots, r$. 
		$T_j$ are the $(n + \nu_j) \times (n + \nu_{j - 1})$ truncations of $U_j$ for $j = 1,2, \dotsc, r$.
		\end{assumption}
By applying 
Proposition \ref{TruncationTransformation} repeatedly, we see that the squared singular value density of $T_r T_{r - 1} \dotsm T_1 X$ is a polynomial ensemble, as shown in \cite{Kieburg-Kuijlaars-Stivigny15}. 
Then we have the following result as another consequence of Lemma \ref{Transformation-phi}.

\begin{corollary} \label{cor:truncated_iterated}
  Let $T_1, \dotsc, T_r$ be the truncated unitary matrices from Assumption \ref{defn:T_matrices} and $X$ be the random matrix 
	from Assumption  \ref{defn:X_matrix} with $\ell=n+\nu_0$, independent of $T_1, \dotsc, T_r$. Then the squared singular value density of the product $T_r \dotsm T_1 X$ is a polynomial ensemble on $[0, \infty)$ with the biorthogonal system consisting of polynomials
  \begin{equation} \label{TruncationPk}
    P_k(x) = \frac{1}{2\pi i} \oint_{\Sigma} 
    \MeijerG{1,r+1}{r+1,r+1}{0,-\nu_1-\mu_1, \ldots, -\nu_r-\mu_r}{0,-\nu_1, \ldots, -\nu_r}{-s} p_k \left( \frac{x}{s} \right) \frac{ds}{s}
  \end{equation}
  where $\Sigma$ is a simple counter-clockwise contour around $0$ lying inside the unit disk, and dual functions
  \begin{equation} \label{TruncationQk} 
    Q_k(y) = \int_0^1  \MeijerG{r,0}{r,r}{\nu_1+\mu_1, \ldots, \nu_r + \mu_r}{\nu_1, \ldots, \nu_r}{t}
    q_k \left( \frac{y}{t} \right) \frac{dt}{t}.
  \end{equation}
  Furthermore, the correlation kernel for the transformed polynomial ensemble is
  \begin{multline} \label{trancationK_n}
    K_n(x,y) = \frac{1}{2\pi i} \oint_{\Sigma} \frac{ds}{s} \int_0^{1} \frac{dt}{t}
    \MeijerG{1,r+1}{r+1,r+1}{0,-\nu_1-\mu_1, \ldots, -\nu_r-\mu_r}{0, -\nu_1, \ldots, -\nu_r}{-s} \\
    \times \MeijerG{r,0}{r,r}{\nu_1+\mu_1, \ldots, \nu_r+\mu_r}{\nu_1, \ldots, \nu_r}{t}
    k_n \left(\frac{x}{s}, \frac{y}{t} \right).
  \end{multline}
\end{corollary}
We will prove Corollary \ref{cor:truncated_iterated} in Section \ref{section: proof corollary trunc}.

We remark that Corollaries \ref{cor:Ginibre_iterated} and \ref{cor:truncated_iterated} can  
also be proved by repeated application of Theorems \ref{GinibreTransformationKernel} and 
\ref{TruncationTransformationKernel}.
However, this would lead to cumbersome manipulations of integrals with Meijer G-functions,
which is avoided by the use of Lemma \ref{Transformation-phi}.

We can use Corollary \ref{cor:truncated_iterated} to study the squared singular values
of a product of the form $T_r \cdots T_1$.
To do so, we can proceed in two ways.
We can either start with a truncated unitary matrix $T_1$
whose squared singular values are a determinantal point process and then apply
Corollary \ref{cor:truncated_iterated} with $X = T_1$, which is multiplied by $r-1$
truncated unitary matrices $T_{r} \cdots T_2$. In this way we recover the result of 
\cite[Corollary 2.6 and Proposition 2.7]{Kieburg-Kuijlaars-Stivigny15}, as we discuss briefly in 
Remark \ref{eq:truncated_KKS}. (In Section \ref{subsec:product_Ginibre_KZ} we consider an analogous problem, the squared singular values of $G_r \dotsm G_1$, in more detail and recover the result of \cite[Proposition 5.1]{Kuijlaars-Zhang14}.)

The squared singular values of $T_1$ are a determinantal point process only if $\mu_1 \geq n$, 
in which case it is the Jacobi ensemble \eqref{Jacobi-ensemble}
with exponents $\nu_1$ and $m_1 - 2n - \nu_1 = \mu_1 - n \geq 0$.
If $\mu_1 < n$, then $T_1$ has a singular value at $1$ of multiplicity $\geq n-\mu_1$,
and the squared singular values are not a determinantal process in the usual, non-degenerate sense. 

However, we may use Corollary \ref{cor:truncated_iterated} in a second way, by 
first taking a limit where the random matrix $X$ approaches an $n \times n$ deterministic matrix 
with distinct singular values, say
\begin{equation} \label{eq:defn_A}
  A = \diag(\sqrt{a_1}, \sqrt{a_2}, \dotsc, \sqrt{a_n}), \quad a_j \in (0, \infty) \text{ are distinct}. 
\end{equation}
Then the squared singular values of $T_r \dotsm T_1 A$ are a degenerate form of a polynomial ensemble,
and we discuss this in Section \ref{subsec:proof_degenerate_Ginibre}.
in the context of products with complex Ginibre matrices.

In a further step we can take the limit $A \to I$, which in the context of Corollary \ref{cor:truncated_iterated}
leads to the final main result of this paper. It  answers a question posed in  \cite[Section 7]{Kieburg-Kuijlaars-Stivigny15}, 
about the precise conditions on truncated unitary matrices so that the squared singular values of 
$ T_r \cdots T_1$ are a determinantal point process. 

\begin{theorem} \label{TruncationProduct} 
  Let $T_1, \dotsc, T_r$ be defined as in Assumption \ref{defn:T_matrices} with $\nu_0 = 0$. 
	Then the squared singular values of $T_r \cdots T_1$ are a determinantal point process if and only if 
  \begin{equation} \label{TruncationProductCondition} 
    n \leq \sum_{j=1}^r \mu_j. 
  \end{equation}
	
	If \eqref{TruncationProductCondition} does not hold, then $T_r \cdots T_1$ has a singular value at $1$
	of multiplicity $\geq n - \sum_{j=1}^r \mu_j$.
  
  If \eqref{TruncationProductCondition} holds, then the determinantal point process is a polynomial ensemble
	on $[0,1]$ with biorthogonal system consisting of polynomials
  \begin{equation} \label{Pk-ProductTruncation}
    \begin{split}
      P_k(x) & = \sum_{j=0}^k (-1)^{k-j} \binom{k}{j} \prod_{l=1}^r \frac{(j+\nu_l+\mu_l)!}{(j+\nu_l)!}  x^j \\
      & = k! \MeijerG{0,r+1}{r+1,r+1}{k+1,-\nu_1-\mu_1, \ldots, -\nu_r-\mu_r}{0,-\nu_1, \ldots, -\nu_r}{x}
    \end{split}
  \end{equation}
  and dual functions
  \begin{equation} \label{Qk-ProductTruncation}
    Q_k(y) =  \frac{1}{k!} \MeijerG{r+1,0}{r+1,r+1}{-k,\nu_1+\mu_1, \ldots, \nu_r + \mu_r}{0,\nu_1, \ldots,\nu_r}{y},
  \end{equation}
  and the correlation kernel is given by
  \begin{equation} \label{Kn-ProductTruncation}
    \begin{split}
      K_n(x,y)  = {}& \frac{1}{(2\pi i)^2} \int_C ds \int_{\gamma} dt
      \prod_{j=0}^r \frac{\Gamma(s+1+\nu_j) \Gamma(t+1+\nu_j+\mu_j)}{\Gamma(t+1+\nu)\Gamma(s+1+\nu_j+\mu_j)}
      \frac{x^t y^{-s-1}}{s-t} \\ 
      = {}& - \int_0^1 \MeijerG{0,r+1}{r+1,r+1}{-\nu_0-\mu_0, \ldots, -\nu_r-\mu_r}{-\nu_0, \ldots, \nu_r}{ux} \\
      & \phantom{- \int_0^1} \times
      \MeijerG{r+1,0}{r+1,r+1}{\nu_0+\mu_0,\ldots,\nu_r+\mu_r}{\nu_0, \ldots, \nu_r}{uy} du
    \end{split}
  \end{equation}
  where $\nu_0 = 0$ and $\mu_0 = -n$. The contour $C$ is a positively oriented Hankel contour in the
  left-half of the complex $s$-plane that starts and ends at $-\infty$ and encircles $(-\infty,-1]$, and 
  $\gamma$ is a closed contour around $[0,n]$ that is disjoint from $C$.
\end{theorem}

We will prove Theorem \ref{TruncationProduct} in Section \ref{section: proof TruncationProduct}.

\begin{rmk}
 In a recent paper \cite{Imamura-Sasamoto15}, Imamura and Sasamoto derived a polynomial ensemble related to the O'Connell--Yor directed random polymer model. The biorthogonal system of the polynomial ensemble, see \cite[Formulas (2.1) and (2.3)]{Imamura-Sasamoto15}, has similarities with our formulas \eqref{Pk} and \eqref{Qk}.
\end{rmk}
\begin{rmk}
  In this paper we consider the product $G_r G_{r - 1} \dotsm G_1 X$ or $T_r T_{r - 1} \dotsm T_1 X$ with the parameter $r$ fixed. It is possible to consider the correlation of the squared singular values of all the products $G_r G_{r - 1} \dotsm G_1 X$ or $T_r T_{r - 1} \dotsm T_1 X$ with $r = 0, 1, 2, \dotsc$, see the recent paper \cite{Strahov15} by Strahov.
\end{rmk}

\subsection*{Outline}

In Section \ref{section: GUE}, we prove Theorem \ref{theorem GUE jpdf} and Theorem \ref{theorem: kernel GUE} about the eigenvalues of the sum of a GUE matrix with another Hermitian random matrix. We will also discuss some concrete consequences of those results.
In Section \ref{sec:general_result}, we prove Lemma \ref{Transformation-phi}.
In Section \ref{section: Ginibre}, we apply Lemma \ref{Transformation-phi} to prove the results in Theorem \ref{GinibreTransformationKernel} on multiplication with a complex Ginibre matrix and corollaries on multiplication with several Ginibre matrices. 
In Section \ref{section: trunc}, we apply Lemma \ref{Transformation-phi} to prove the results on multiplication with a truncated unitary matrices stated in Theorem \ref{TruncationTransformationKernel}, and the corollaries on multiplication with several truncated matrices. At last we prove Theorem \ref{TruncationProduct}. 

\section{Addition of a GUE matrix} \label{section: GUE}

\subsection{Proof of Theorem \ref{theorem GUE jpdf}} \label{section: proof GUE jpdf}

Let $H$ be an $n \times n$ GUE matrix, and temporarily let $M$ be a fixed $n \times n$ Hermitian matrix with eigenvalues $x_1, \ldots, x_n$. Then the eigenvalue density of $H$ is given by \eqref{eq:distr_GUE} and the random matrix $Y=H+M$ has distribution
\begin{equation}
  \frac{1}{Z_n}e^{-\frac{1}{2}\Tr\left(Y-M\right)^2}dY, \qquad Z_n = (2\pi)^{\frac{n}{2}} \prod^n_{k = 1} k!,
\end{equation}
which is the distribution for a random matrix $Y$ in the GUE with external source $M$ \cite{Brezin-Hikami96}.
It is known from \cite{Bleher-Kuijlaars05}, \cite{Brezin-Hikami96} and \cite{Johansson01a}
that the eigenvalues $y_1,\ldots, y_n$ of $Y=H+M$ are distributed according   
\begin{equation} \label{eq:GUE_ext_jpdf}
  P(y; x) = \frac{1}{\widehat Z_n\Delta_n(x)} \Delta_n(y) \det \left[ e^{-\frac{1}{2} (y_j - x_k)^2} \right]^n_{j, k = 1}, \qquad \widehat{Z}_n = n! (2\pi)^{\frac{n}{2}},
\end{equation} 
where $x=(x_1,\ldots, x_n)$ and $y=(y_1,\ldots, y_n)$.

Now, let $M$ be the random matrix satisfying Assumption \ref{defn:M_matrix}. We obtain the joint probability density function of the eigenvalues $y_1,\ldots, y_n$ of $H+M$ by integrating \eqref{eq:GUE_ext_jpdf} over the eigenvalues $x_1,\ldots, x_n$ of $M$. If their density is the polynomial ensemble \eqref{f-pol-ensemble}, we obtain that the density of the eigenvalues $y_1, \dotsc, y_n$ of $H + M$ is ($Z_n$ is defined in \eqref{f-pol-ensemble}, the polynomial ensemble of the eigenvalues of $M$)
\begin{equation} 
  \begin{split}
    P(y) = {}& \int_{\mathbb R^n} P(y;x) \Delta_n(x)\det \left[ f_{k-1}(x_j) \right]_{j,k=1}^n  dx_1 \dotsm dx_n \\
    = {}& \frac{\Delta_n(y)}{Z_n\widehat Z_n }\int_{\mathbb R^n} \det \left[ f_{k - 1}(x_j) \right]^n_{j, k = 1} 
    \det \left[ e^{-\frac{1}{2}(y_j - x_k)^2} \right]^n_{j, k = 1}  dx_1 \dotsm dx_n.
  \end{split}
\end{equation}
By the \Andreief\ formula, see e.g.\ \cite{Deift-Gioev09}, we have
\begin{equation} \label{eq:jpdf_calculated}
  P(y) = \frac{1}{Z_n' } \Delta_n(y)\det \left( \int_{-\infty}^{\infty} f_{k - 1}(x) e^{-\frac{1}{2}(y_j - x)^2} dx \right)^n_{j, k = 1},
\end{equation}
with $Z_n' = n! Z_n \widehat{Z}_n$, which is the polynomial ensemble defined by \eqref{f-pol-ensemble 2} and \eqref{def F GUE}.
\qed

\subsection{Proof of Theorem \ref{theorem: kernel GUE}} \label{section: proof GUE}

We assume that $p_k$ and $q_k$ for $k=0,1,\ldots, n-1$, are the biorthogonal system 
as in Assumption \ref{defn:M_matrix} associated to $M$. In other words, we assume that for 
$k=0,\ldots, n-1$, $p_k$ is a polynomial of degree $k$, $q_k$ is in the linear span of 
$f_0,\ldots, f_{n-1}$, and they satisfy \eqref{eq:pqbiort}.

Denote the Weierstrass transform \cite{Hirschman-Widder55} of a function $\varphi$ by
\begin{equation}\label{def Weierstrass}
\mathcal W\varphi(y)=\frac{1}{\sqrt{2\pi}}\int_{-\infty}^{\infty}\varphi(t)e^{-\frac{1}{2}(y-t)^2}dt.
\end{equation}
The integral transform
\begin{equation}\label{def Weierstrass inv}
\mathcal W^{-1}\Phi(x)=\frac{1}{\sqrt{2\pi} i}\int_{-i\infty}^{+i\infty}\Phi(s)e^{\frac{1}{2}(x-s)^2}ds
\end{equation}
is the inverse Weierstrass transform: we have $\left(\mathcal W\circ \mathcal W^{-1}\right)\Phi=\Phi$ for a large class 
of functions. We will need this inversion property only for polynomials $\Phi$. A property of the Weierstrass transform 
is that it sends monic polynomials to monic polynomials of the same degree.

Formulas \eqref{Pk} and \eqref{Qk} can now be written as
\begin{equation} \label{PkandQk-Weierstrass}
P_k(x)=\mathcal W^{-1}p_k(x),\qquad Q_k(y)=\mathcal W q_k(y).
\end{equation}
Since $q_k$ is in the linear span of $f_0,\ldots, f_{n-1}$, $Q_k$ is in the linear span of  $\mathcal W f_0,\ldots, \mathcal W f_{n-1}$, 
and thus also in the linear span of $F_0,\ldots, F_{n-1}$ since $F_j=\sqrt{2\pi}\mathcal W f_j$ by \eqref{def F GUE}.

By \eqref{PkandQk-Weierstrass}, \eqref{def Weierstrass} and Fubini's theorem, we obtain
\begin{equation} 
  \begin{split}
    \int_{-\infty}^{\infty} P_j(x) Q_k(x) dx  = {}&
    \int_{-\infty}^{\infty}\mathcal W^{-1}p_j(x)\, \mathcal W q_k(x)dx\\
    = {}& \int_{-\infty}^{\infty}   \left(\frac{1}{\sqrt{2\pi}} \int_{-\infty}^{\infty} \mathcal W^{-1} p_j(x)  
      e^{-\frac{1}{2}(t-x)^2} dx  \right) q_k(t) dt.
  \end{split}
\end{equation}
The expression within parentheses is the Weierstrass transform of $\mathcal W^{-1}p_j$ and is thus equal to $p_j(t)$.
Hence 
\begin{equation}\label{orth Weierstrass}
	\int_{-\infty}^{\infty} P_j(x) Q_k(x) dx = \int_{-\infty}^{\infty} p_j(x)q_k(x)dx= \delta_{j,k} 
\end{equation}
by \eqref{eq:pqbiort}. This proves part \ref{enu:theorem: kernel GUE_a} of the theorem.

\medskip
Next, we substitute \eqref{Pk} and \eqref{Qk} into \eqref{eq:Kn hat} and use \eqref{eq:Kn}. 
This gives \eqref{integral transform} and proves part \ref{enu:theorem: kernel GUE_b}.
\medskip

Formula \eqref{orth Weierstrass} is also valid for $j=n$ and $k = 0, 1, \ldots,n-1$.
Thus $P_n$ is a monic polynomial of degree $n$ which is orthogonal to $Q_0, \ldots, Q_{n-1}$.
This implies that $P_n$ is the average characteristic polynomial in the ensemble \eqref{f-pol-ensemble 2}
and we have proven part \ref{enu:theorem: kernel GUE_c} of the theorem. \qed


\subsection{Example: GUE plus a unitary invariant random matrix} \label{subsubsec:GUE+1MM}

In this subsection we assume that the random matrix $M$ from Assumption \ref{defn:M_matrix} is defined by the 
unitary invariant probability measure \cite{Deift99}
  \begin{equation}\label{UE}
    \frac{1}{C_n} e^{-\Tr V(M)}dM,\qquad C_n=\int e^{-\Tr V(M)}dM,
  \end{equation}
  where $V$ is a real-valued function such that the integral defining $C_n$ is convergent 
	and which is allowed to depend on $n$. This is a random matrix ensemble which is invariant 
	under unitary conjugation. If $V(x) = x^2/2$, it is simply the GUE.
	The joint probability 
	density function for the eigenvalues of $M$ is the polynomial ensemble with $f_k$ given by 
	$f_k(x)=x^ke^{-V(x)}$ in \eqref{f-pol-ensemble}.
  For the biorthogonal systems, it is convenient to let $p_k$ be the monic degree $k$ orthogonal polynomial 
	with respect to the weight $e^{-V(x)}$ on the real line and 
  \begin{equation}
    q_k(x)=\frac{1}{h_k}p_k(x)e^{-V(x)},\qquad h_k=\int_{-\infty}^{\infty}p_k(x)^2e^{-V(x)}dx.
  \end{equation} 
The correlation kernel $k_n$ can be written as
\begin{equation}\label{1MMkernelsum1}
  k_n(x,y) = e^{-V(y)}\sum_{j=0}^{n-1}\frac{1}{h_j}p_j(x)p_{j}(y),
\end{equation}
and by the Christoffel-Darboux formula we also have
\begin{equation}\label{1MMkernel}
k_n(x,y)=\frac{1}{h_{n-1}}e^{-V(y)}\frac{p_n(x)p_{n-1}(y)-p_n(y)p_{n-1}(x)}{x-y}.
\end{equation}

Theorem \ref{theorem: kernel GUE} applies to this case and we obtain the formula 
\begin{multline}\label{kernelGUE+UE}
K_n(x,y)=\frac{1}{2\pi ih_{n-1}} \int_{-i\infty}^{+i\infty}ds\int_{-\infty}^{\infty}dt\, 
	\frac{p_n(s)p_{n-1}(t)-p_n(t)p_{n-1}(s)}{s-t}\\
	\times\quad  e^{-V(t)} e^{\frac{1}{2}((s-x)^2-(t-y)^2)},
\end{multline}
for the eigenvalue correlation kernel of the matrix $Y=H+M$, where $H$ is a GUE matrix of size $n$ and $M$ has the distribution \eqref{UE}.
We show how to rewrite this formula in a way that may be suitable for asymptotic analysis.

For a large class of potentials $V$ (in particular, for $n$-dependent potentials of the form $V(x)=nV_0(x)$ with 
$V_0$ independent of $n$), large $n$ asymptotics for $p_n(z)$ and $p_{n-1}(z)$ are known for $z$ anywhere in the 
complex plane. However, it is not straightforward to apply saddle point techniques on the integrals in 
\eqref{kernelGUE+UE}, especially because one has to integrate over the real line, where the zeros of the orthogonal polynomials 
$p_n$, $p_{n-1}$ are, and where the integrand is oscillatory. Therefore, we derive an alternative expression 
for the correlation kernel, which involves Cauchy transforms of the orthogonal polynomials and which avoids 
integration over the real line.

Define a $2\times 2$ matrix-valued function (which is the solution to the Riemann-Hilbert problem for orthogonal polynomials
\cite{Deift99}, \cite{Fokas-Its-Kitaev92})
\begin{equation}
Y(z)=\begin{pmatrix}p_n(z)&\frac{1}{2\pi i}\int_{-\infty}^{\infty}p_n(s)e^{-V(s)}\frac{ds}{s-z}\\
-2\pi ih_{n-1}^{-1}p_{n-1}(z)&-h_{n-1}^{-1}\int_{-\infty}^{\infty}p_{n-1}(s)e^{-V(s)}\frac{ds}{s-z}
\end{pmatrix}.
\end{equation}
The entries in the first column are defined everywhere in the complex plane, the Cauchy transforms 
in the second column are defined for  $z\in\mathbb C\setminus \mathbb R$.
We can express the correlation kernel $k_n$ in \eqref{1MMkernel} in terms of $Y$. 
It is straightforward to check that \cite[Chapter 8]{Deift99}
\begin{equation}
	k_n(s,t)=\frac{-1}{2\pi i(s-t)}e^{-V(t)}\begin{pmatrix}0&1\end{pmatrix}Y^{-1}(s)Y(t)\begin{pmatrix}1\\0\end{pmatrix},
\end{equation}
for $s\in\mathbb C$, $t\in\mathbb R$.
Writing $Y_+(t)$ (resp.\ $Y_-(t)$) for the limit of $Y(z)$ as $z$ approaches $t\in\mathbb R$ from 
the upper half plane (resp.\ lower half plane), we have the relation 
\begin{equation}
Y_+(t)=Y_-(t)\begin{pmatrix}1&e^{-V(t)}\\0&1\end{pmatrix},\qquad t\in\mathbb R.
\end{equation}
It follows that 
\begin{equation}
e^{-V(t)}Y(t)\begin{pmatrix}1\\0\end{pmatrix}=
\left(Y_+(t)-Y_-(t)\right)\begin{pmatrix}0\\1\end{pmatrix},\qquad t\in\mathbb R,
\end{equation}
and that
\begin{equation}
k_n(s,t)=\frac{-1}{2\pi i(s-t)}\begin{pmatrix}0&1\end{pmatrix}Y^{-1}(s)\left(Y_+(t)-Y_-(t)\right)\begin{pmatrix}0\\1\end{pmatrix},
\end{equation}
for $t\in\mathbb R$ and for any $s\in\mathbb C$.
In \eqref{integral transform}, we need to integrate $k_n(s,t)$ in $t$ over the real line. This integral
can be deformed to a contour $\Gamma$ which consists of a curve in the lower half plane 
oriented from left to right (for example $\Gamma = \mathbb R - i \delta$ for some $\delta >0$), 
and its complex conjugate in the upper half plane oriented from right to left.
We obtain, for any $s$ which is not on $\Gamma$,
\begin{multline}
\int_{-\infty}^{\infty}  k_n(s,t) e^{-\frac{1}{2}(y-t)^2} dt 
	= \frac{1}{2\pi i} \int_{\Gamma} \frac{1}{s-\zeta} 
		\begin{pmatrix} 0& 1 \end{pmatrix} Y^{-1}(s) Y(\zeta) \begin{pmatrix} 0 \\ 1 \end{pmatrix} e^{-\frac{1}{2}(y- \zeta)^2}  d \zeta \\
			+ \begin{cases}		
				e^{-\frac{1}{2}(y-s)^2} & \text{if $s$ is inside $\Gamma$,} \\
				0 & \text{if $s$ is outside $\Gamma$.}
				\end{cases}
\end{multline}
The last term is the residue contribution from the pole at $\zeta = s$. In \eqref{integral transform}, the integration over the imaginary axis can be changed to any vertical line $C$ as contour
of integration for $s$. Let $w_{\pm}$ be the intersection points of $C$
with $\Gamma$, with $w_+$ in the upper half plane and $w_- = \overline{w}_+$. (We assume there are only two intersection points.)

Then the contribution from the additional term to ${K}_n(x,y)$  is
\begin{equation} 
  \begin{split}
    &\frac{1}{2\pi i} \int_{w_-}^{w_+} 	e^{-\frac{1}{2}(y-s)^2}
    e^{\frac{1}{2}(x-s)^2} ds \\
    &\qquad = \frac{1}{2\pi i (y-x)} e^{\frac{1}{2}(x^2 - y^2)} 
    \left(
      e^{(y- x)w_+} - e^{(y- x)w_-} \right) \\
    &\qquad 	= \frac{1}{\pi (y-x)} e^{\frac{1}{2}(x^2 - y^2)} 
    e^{(y- x) \Re w_+} \sin \left( (y- x) \Im w_+\right),
  \end{split}
\end{equation}
where $x$ and $y$ are real.

We thus obtain the alternative expression for the eigenvalue correlation kernel $K_n$ of $H+M$,
\begin{multline}\label{kernelGUE+UECauchy}
{K}_n(x,y) = 
	\frac{1}{\pi (y-x)} e^{\frac{1}{2}(x^2 - y^2)} 
		e^{(y- x) \Re w_+} \sin \left( (y- x) \Im w_+\right)  \\
				+
	\frac{1}{(2\pi i)^2} \int_C ds \int_{\Gamma} \frac{d\zeta}{s-\zeta} 
		\begin{pmatrix} 0 & 1 \end{pmatrix} Y^{-1}(s) Y(\zeta) \begin{pmatrix} 0 \\ 1 \end{pmatrix} e^{\frac{1}{2}((x-s)^2-(y-\zeta)^2)},
\end{multline}
in which integration over the real line is avoided. This formula may be more convenient for 
asymptotic analysis than \eqref{kernelGUE+UE}, but we do not aim to investigate this here.

\begin{rmk}\label{rmk 2MM}
The random matrix $H+M$ appears in the two-matrix model, which is defined as a 
measure on the space of pairs $(M_1,M_2)$ of Hermitian $n\times n$ matrices, 
\begin{equation}\label{2MM}
\frac{1}{C_n}e^{-\Tr\left(W_1(M_1)+W_2(M_2)- M_1M_2\right)}dM_1dM_2,
\end{equation}
for certain functions $W_1$ and $W_2$ such that the above defines a probability distribution.
If $W_1(M_1) = M_1^2/2$, as observed in \cite{Duits14}, the probability measure \eqref{2MM} can be written as
\begin{equation}\label{2MMb}
\frac{1}{C_n}e^{-\Tr\left(\frac{1}{2}(M_1-M_2)^2+W_2(M_2)-\frac{1}{2}M_2^2\right)}d(M_1-M_2)dM_2.
\end{equation}
This implies that $M_2$ and $M_1- M_2$ are independent random matrices: $M_1-M_2$ is a GUE matrix and $M_2$ is a random matrix from a unitary invariant one-matrix model with probability distribution (\ref{UE}), with 
 \begin{equation}\label{V}V(M)=W_2(M)-\frac{1}{2}M^2.
 \end{equation} 
In other words, the matrix $M_1$ in the two-matrix model \eqref{2MM} then takes the form $M_1=H+M$, where $H$ is a GUE matrix and $M$ is a random matrix from the ensemble \eqref{UE}.
\end{rmk}

\begin{rmk} 
If $V$ in \eqref{UE} is a polynomial, the eigenvalue correlation kernel $K_n$ for $H+M$ can be expressed in terms of multiple orthogonal polynomials, since we can interpret $H+M$ as the matrix $M_1$ in the two-matrix model, but the asymptotic analysis of those multiple orthogonal polynomials using Riemann-Hilbert techniques is hard in general, in particular if the degree of $V$ is large.
The case where $V$ is a quartic symmetric polynomial has been studied in detail in \cite{Duits-Geudens13}.
Our results yield an alternative expression for the correlation kernel of the eigenvalues of $H+M$, for a general potential $V$, which does not involve multiple orthogonal polynomials but only usual orthogonal polynomials and contour integrals thereof. 
 \end{rmk}

\section{Proof of Lemma \ref{Transformation-phi}} \label{sec:general_result}

\paragraph{Part \ref{enu:ransformation-phi_a}}

From the definitions \eqref{Fkphi-definition} and \eqref{eq:Qkphi} it is clear that $Q_k$ belongs to
the linear span of $F_0, \ldots, F_{n-1}$, since $q_k$ belongs to the linear
span of $f_0, \ldots, f_{n-1}$. It is also clear from \eqref{eq:Pkphiseries} 
that $P_k$  is a polynomial of degree $k$ for every $k=0, \ldots, n-1$.

To verify the biorthogonality
we define an operator $\mathcal L$ acting on polynomials $p$ by
\begin{equation} \label{def:L} 
	\mathcal L p(x)  = \sum_{j=0}^{\deg(p)} a_j b_j x^j \qquad 
	\text{where} \quad p(x) = \sum_{j=0}^{\deg(p)} a_j x^j 
	\end{equation}
is arbitrary while $b_j$ are defined in \eqref{eq:bj}, and an operator $\mathcal M$ acting on functions $q$ on $[0,\infty)$ by
\begin{equation} \label{def:M} 
	\mathcal M q(y) = \int_0^{\infty} \varphi(t) q \left( \frac{y}{t} \right) \frac{dt}{t}, \qquad y > 0. 
	\end{equation}
Then $P_k = \mathcal L p_k$ and $Q_k = \mathcal M q_k$ by  \eqref{eq:Pkphiseries} and \eqref{eq:Qkphi}.
We prove the identity
\begin{equation} \label{LM identity}  
	\int_{0}^{\infty} \mathcal L p (x)  \,  \mathcal M q(x) \, dx =
	\int_0^{\infty} p(x) q(x) \, dx
	\end{equation}
and then the biorthogonality follows from 
\begin{equation}
  \int_0^{\infty} P_j(x) Q_k(x) \, dx = \int_0^{\infty} p_j(x) q_k(x) \, dx = \delta_{j,k}
\end{equation}
since $P_j = \mathcal L p_j$, $Q_k = \mathcal M q_k$ for $j,k =0, 1, \ldots, n-1$.

To prove \eqref{LM identity}  we calculate by Fubini's theorem for $k=0, 1, \ldots$
\begin{equation}
  \int_{0}^{\infty} x^k  \,  \mathcal M q(x) \, dx
  = \int_0^{\infty} \varphi(t) \left( \int_0^{\infty}  x^k 
    q \left(\frac{x}{t}\right) dx \right) \frac{dt}{t}.
\end{equation}
We substitute $x = tu$ in the inner integral and obtain
\begin{align} \nonumber 	
	\int_{0}^{\infty} x^k  \,  \mathcal M q(x) \, dx
	& =  \int_0^{\infty} t^{k} \varphi(t) dt \, 
		\int_0^{\infty} u^k  q(u) du  \\
		& = b_k^{-1}  	\int_0^{\infty} x^k  q(x) dx. \label{xkM integral}
		\end{align}
By the definition \eqref{def:L} we have $\mathcal L x^k = b_k x^k$, and thus  \eqref{xkM integral} implies that \eqref{LM identity} holds 
if $p(x) = x^k$ for $k=0,1, \dotsc$. By linearity \eqref{LM identity} holds for every polynomial $p$. 

\paragraph{Part \ref{enu:ransformation-phi_b}}

The transformation \eqref{eq:Pkphiseries} is the Hadamard (or termwise)
product of  $p_k$ with the function $\psi$. 
The Hadamard product of two convergent Laurent series
$a(x) = \sum_{j=-\infty}^{\infty} a_j x^j$ and $b(x) = \sum_{j=-\infty}^{\infty} b_j x^j$
has a well-known (and easy to prove) contour integral  representation, namely 
\begin{equation} \label{Hadamard Integral} 
	\sum_{j=-\infty}^{\infty} a_j b_j x^j 
	 = \frac{1}{2\pi i} \oint_{\Sigma_r} a(s) b \left( \frac{x}{s} \right) \frac{ds}{s},
	\end{equation}
where $\Sigma_r$ is the circle of radius $r$ around the origin with positive
orientation. Here we assume that the Laurent series for $a(x)$ converges
for $r_a < |x| < R_a$, and the  Laurent series for $b(x)$ converges for $r_b < |x| < R_b$. Then formula \eqref{Hadamard Integral} is 
valid for $r_a r_b < |x| < R_a R_b$ and $r \in (r_a, R_a) \cap 
(\frac{|x|}{R_b}, \frac{|x|}{r_b})$.

The formula \eqref{eq:Pkphiintegral} follows from  \eqref{Hadamard Integral}
and the definitions \eqref{eq:Pkphiseries} and \eqref{eq:psi}.

\paragraph{Part \ref{enu:ransformation-phi_c}}

The correlation kernel for \eqref{Fphi-pol-ensemble} is
\begin{equation}
  K_n(x,y) = \sum_{k=0}^{n-1} P_k(x) Q_k(y).
\end{equation}
We insert the integral representations \eqref{eq:Pkphiintegral} and \eqref{eq:Qkphi} for $P_k$ and $Q_k$,
and interchange the sum with the integrals. Then \eqref{eq:Knhatphi} follows because
of \eqref{eq:Kn}. \qed

\section{Multiplication with complex Ginibre matrices}\label{section: Ginibre}

\subsection{Proof of Theorem \ref{GinibreTransformationKernel}}\label{section: proof Ginibre}

Theorem \ref{GinibreTransformationKernel} follows immediately from Lemma \ref{Transformation-phi}. 
In the situation of Proposition \ref{GinibreTransformation} we have
\begin{equation}
  \varphi(t) = t^{\nu} e^{-t}
\end{equation}
with $\nu$ a non-negative integer. Then the moments \eqref{eq:bj} are
\begin{equation}
  b_j = \left[ \int_0^{\infty} t^{j+\nu} e^{-t} dt \right]^{-1}
  = \left[ (j+\nu)! \right]^{-1},
\end{equation}
and for  the function $\psi$ from \eqref{eq:psi} we take
\begin{equation}
  \psi(x) = \sum_{j=-\nu}^{\infty} \frac{x^j}{(j+\nu)!} = x^{-\nu} e^x.
\end{equation}
Then parts \ref{enu:ransformation-phi_a} and \ref{enu:ransformation-phi_b} of Lemma \ref{Transformation-phi} give the formulas \eqref{eq:biorthogonal_pd_P} and \eqref{eq:biorthogonal_pd_Q} for the biorthogonal system in Theorem \ref{GinibreTransformationKernel}\ref{enu:GinibreTransformationKernel_a}, and Theorem \ref{GinibreTransformationKernel}\ref{enu:GinibreTransformationKernel_c} is also derived in the same way. Lemma \ref{Transformation-phi}\ref{enu:ransformation-phi_c} gives the transformed kernel \eqref{eq:Knhatphi} in Theorem \ref{GinibreTransformationKernel}\ref{enu:GinibreTransformationKernel_b}. \qed

\subsection{Proof of Corollary \ref{cor:Ginibre_iterated}}
\label{section: corollary Ginibre}
Assume $X$ and $G_1, \dotsc, G_r$ are as in Corollary \ref{cor:Ginibre_iterated}. 
Applying Proposition \ref{GinibreTransformation} $r$ times, see also \cite{Kuijlaars-Stivigny14}, 
we find that the squared singular value density of $Y = G_r \cdots G_1 X$ is a polynomial ensemble with functions  
$F_k$ that are the Mellin convolution of $f_k$ with
\begin{equation} \label{eq:varphi_multiple_G}
  \varphi = \varphi_r \ast \varphi_{r-1} \ast \cdots \ast \varphi_1
\end{equation}
where $\varphi_j(t) = t^{\nu_j} e^{-t}$ for $j = 1, \ldots, r$, and $\ast$ denotes the Mellin convolution. Thus $\varphi$ is a Meijer G-function
\begin{equation} \label{phi:MeijerG} 
	\varphi(t) = \MeijerG{r,0}{0,r}{-}{\nu_1, \ldots, \nu_r}{t} 
	\end{equation}
which has the moments
\begin{equation} \label{bj:MeijerG} 
	b_j^{-1} = \int_0^{\infty} t^j \varphi(t) dt 	= \prod_{k=1}^r (j + \nu_k)!. 
	\end{equation}

Now we apply Lemma \ref{Transformation-phi}. We take for $\psi$
\begin{equation} \label{psi:defseries} 
	\psi(s) = \sum_{j=0}^{\infty} \frac{s^j}{\prod_{k=1}^r (j+\nu_k)!}, 
	\end{equation}
which is a generalized hypergeometric function
\begin{equation}
  \psi(s) = \frac{1}{\prod_{k=1}^r \nu_k!} 
  \hyperF{1}{r}{1}{\nu_1 + 1, \dotsc, \nu_r + 1}{s}
\end{equation}
and it is also a Meijer G-function
\begin{equation} \label{psi:MeijerG} 
	\psi(s) = \MeijerG{1,1}{1,r+1}{0}{0, -\nu_1, \ldots, -\nu_r}{-s}. 
	\end{equation}
Then Corollary \ref{cor:Ginibre_iterated} is a direct consequence of Lemma \ref{Transformation-phi}. \qed

\subsection{Products of Ginibre  matrices} \label{subsec:product_Ginibre_KZ}

We consider the product 
\begin{equation}
  Y_r = G_r G_{r-1} \cdots G_1
\end{equation}
of independent complex Ginibre matrices where $G_j$ are defined in Assumption \ref{defn:G_matrices} with $\nu_0 = 0$. 
It was shown in \cite{Kuijlaars-Zhang14} that the squared singular value density of $Y_r$ is a polynomial ensemble with correlation kernel
\begin{equation} \label{kernel product}
  K_n(x,y)=\frac{1}{(2\pi i)^2}\int_{-\frac{1}{2}+i\mathbb R}dv\oint_{\gamma}du 
  \prod_{j=0}^r \frac{\Gamma(v+\nu_j+1)}{\Gamma(u+\nu_j+1)}\, \frac{\Gamma(u-n+1)}{\Gamma(v-n+1)} \, \frac{x^uy^{-v-1}}{v - u},
\end{equation}
where $\nu_0=0$ and $\gamma$ encircles $0,1,\ldots, n$ once in positive direction and lies to the right 
of $-\frac{1}{2}+i\mathbb R$, and $\Gamma$ is the Euler Gamma function. For $r=1$, this is the correlation 
kernel for the eigenvalues in the complex Wishart Ensemble $G_1^* G_1$ with parameter $\nu_1$.

We show how to obtain \eqref{kernel product} from Theorem \ref{GinibreTransformationKernel} 
by induction on $r$.  

\subsubsection{Base step: Case $r=1$}
The $r = 1$ case of \eqref{kernel product} differs from the well known double contour 
integral formula for the complex Wishart ensemble (\cite[Section 5.8]{Forrester10})
\begin{equation} \label{eq:traditional_Wishart}
  K_n(x, y) = \frac{1}{(2\pi i)^2} \oint_{\Sigma} du \oint_{\Gamma} dv \frac{e^{xu} v^{n + \nu_1} (u - 1)^n}{e^{yv} u^{n + \nu_1} (v - 1)^n} \frac{1}{u - v},
\end{equation}
where $\Sigma$ and $\Gamma$ are disjoint closed contours both oriented counterclockwise such that $\Sigma$ encloses 
$0$ and $\Gamma$ encloses $1$. The $r = 1$ case of \eqref{kernel product} is less well known, but it also appeared in 
the literature as a special case of the \emph{multiple Laguerre minor process} in \cite[Theorem 3(c)]{Adler-van_Moerbeke-Wang11}
with the fixed time $n$ and special choice of $\alpha_k = n + \nu_1 - k$ for $k = 1, \dotsc, n$, see 
also \cite{Forrester-Wang15} for its relation to the Laguerre Muttalib--Borodin model. To be precise, the relation is
\begin{equation} \label{eq:K_n_r=1_equivalence}
  \left. K_n(x, y) \right\rvert_{r = 1} = \left. \frac{y^{\nu_1}}{x^{\nu_1}} K(n, y; n, x) \right\rvert_{\text{$\alpha_k = n + \nu_1 - k$ for $k = 1, \dotsc, n$}},
\end{equation}
where $K$ is the correlation kernel defined in \cite[Theorem 3(c)]{Adler-van_Moerbeke-Wang11}, 
and the factor $y^{\nu_1}/x^{\nu_1}$ is to conjugate the kernel into the $\sum^{n - 1}_{k = 0} p_k(x) q_k(y)$ form. 
The choice of contour in \cite[Theorem 3(c)]{Adler-van_Moerbeke-Wang11} is different from that in \eqref{kernel product}, 
but they are equivalent by the residue theorem. In comparison, the better known formula \eqref{eq:traditional_Wishart} 
is a special case of the \emph{Wishart minor process} stated in \cite[Theorem 3(b)]{Adler-van_Moerbeke-Wang11}, 
which was originally studied in \cite{Borodin-Peche08}, \cite{Dieker-Warren09} and \cite{Forrester-Nagao11}. 
Both of the two correlation kernels are constructed by 
$\sum^{n - 1}_{k = 0} p_k(x) q_k(y)$ with $p_k(x) = L^{(\nu_1)}_k(x)$ and 
$q_k(y) = k! \Gamma(k + \nu_1 + 1)^{-1} L^{(\nu_1)}_k(y) y^{\nu_1} e^{-y}$, by expressing 
$p_k(x)$ and $q_k(y)$ into contour integral forms, and then summing them up via a telescoping trick. 
The difference lies in the fact that different contour integral formulas for $p_k$ and $q_k$ are used.

\subsubsection{Induction step}

Let $r \geq 2$, and assume that \eqref{kernel product} with $r$ replaced by $r-1$ is the correlation kernel for 
the squared singular values of $Y_{r-1} = G_{r-1} \cdots G_1$. Then $Y_r = G_r Y_{r-1}$, and 
using \eqref{kernel product}  (with $r-1$ instead of $r$) for $k_n$ in the formula 
\eqref{integral transform product} with $\nu = \nu_r$, we obtain a quadruple integral
on $-1/2 + i\realR$, $\gamma$, $\Sigma$, and $[0, \infty)$ for the new correlation kernel.
We then first deform the closed contour $\Sigma$ to a Hankel contour $L$
that comes from $-\infty$ in the lower half-plane, loops around the
negative real axis and goes to $-\infty$ in the upper half-plane.
We can then interchange the order of integration. The $s$- and $t$-integrals are evaluated
explicitly as 
\begin{align} \label{GammaIntegral}
	\frac{1}{2\pi i} \int_{L} s^{-\nu_{r}-u-1} e^{s} ds & 
	=  \frac{1}{\Gamma(u+\nu_{r}+1)}, \\  \label{GammaRecipIntegral}
	\int_0^\infty t^{\nu_{r}+v} e^{-t}dt & = \Gamma(v+\nu_{r}+1),
\end{align}
by the  integral representations for the Gamma function and its reciprocal.
Then the result is a double contour integral on $-1/2 + i\realR$ and $\gamma$, that is precisely \eqref{kernel product} with parameter $r$.

\subsection{Limiting case of Corollary \ref{cor:Ginibre_iterated}} \label{subsec:proof_degenerate_Ginibre}

We may also approach the calculation of the squared singular values of $Y_r = G_r \cdots G_1$
as a limiting case of Corollary \ref{cor:Ginibre_iterated}. 

In Corollary \ref{cor:Ginibre_iterated} we assume $\nu_0 = 0$ and let the random matrix  $X$ approach a fixed $n \times n$ 
matrix with distinct squared singular values $a_1, \dotsc, a_n$, or without loss of generality, 
$X \to A$ as  defined in \eqref{eq:defn_A}. This is a limiting case of a polynomial ensemble with 
functions $f_{k-1}$ that approach the Dirac delta functions $\delta(x-a_k)$ for every $k=1, \ldots, n$.
A limiting biorthogonal system is given by the polynomials
\begin{equation} \label{eq:pkdegenerate}
	p_k(x) = \prod_{j=1}^k (x-a_j), \qquad k = 0, 1, \ldots, n-1 \end{equation}
with $p_0(x) = 1$ and dual functions $q_k$ that are given as a $(k+1) \times (k+1)$ determinant involving Dirac delta functions
\begin{equation} \label{eq:qkdegenerate}
q_k(y) = \frac{1}{\Delta_{k+1}(a)} 
	\begin{vmatrix} 1 & 1 & \cdots & 1 \\
		a_1 & a_2 &  \cdots & a_{k+1} \\
		\vdots & \vdots & \ddots & \vdots \\
		\delta(y-a_1) & \delta(y-a_2) & \cdots & \delta(y-a_{k+1}) \end{vmatrix}.  
		\end{equation}

Alternatively, we can use, instead of $p_k$ and $q_k$, the Lagrange interpolating polynomials
\begin{equation}
  \widetilde{p}_k(x) = \prod_{j\neq k} \frac{x-a_j}{a_k-a_j}
\end{equation}
and the dual functions $\widetilde{q}_k(y) = \delta(y-a_k)$, for $k=1, \ldots, n$. Note that
$\widetilde{p}_k$ has degree $n-1$ for every $k$. The biorthogonality
\begin{equation}
  \int_0^{\infty} p_j(x) q_k(x) dx = \int_0^{\infty} \widetilde{p}_j(x) \widetilde{q}_k(x) dx = \delta_{j,k}
\end{equation}
is easy to verify for both systems.

The correlation kernel for this degenerate ensemble is
\begin{align} \nonumber 
	k_n(x,y) & = \sum_{k=0}^{n-1} p_k(x) q_k(y) = \sum_{k=1}^n \widetilde{p}_k(x) \widetilde{q}_j(y) \\
			& = \left( \prod_{j=1}^n (x-a_j) \right) \sum_{k=1}^n  \left( \prod_{j \neq k} \frac{1}{a_k-a_j} \right)
			\frac{\delta(y-a_k)}{x-a_k}  \end{align}
and with this expression  for $k_n(x,y)$ we calculate \eqref{Knhat-rGinibre}
by first changing variables $s \mapsto x/s$, $t \mapsto y/t$, and then evaluating the $t$-integral
which results in

\begin{multline}
  K_n(x,y) = 
  \frac{1}{2\pi i} \oint_{\Sigma} \frac{ds}{s} 
  \MeijerG{1,1}{1,r+1}{0}{0, -\nu_1, \ldots, -\nu_r}{-\frac{x}{s}} \left( \prod_{j=1}^n (s-a_j) \right) \\
  \times
  \sum_{k=1}^n  \left( \prod_{j \neq k} \frac{1}{a_k-a_j} \right) \frac{1}{a_k(s-a_k)} 
  \MeijerG{r,0}{0,r}{-}{\nu_1, \ldots, \nu_r}{\frac{y}{a_k}}.
\end{multline}

The summation can be written as a contour integral, by the residue theorem, and we obtain the following result.
\begin{cor} \label{cor:degenerate_Ginibre}
  Let $G_1, \dotsc, G_r$ be the Ginibre random matrices defined in Assumption \ref{defn:G_matrices} with $\nu_0 = 0$, and 
	let $A$ be the deterministic matrix defined in \eqref{eq:defn_A}. Then the correlation kernel 
	of the squared singular values of $G_r G_{r - 1} \dotsm G_1 A$ is
  \begin{multline} \label{eq:degenerate_Ginibre}
    K_n(x,y) = 
    \frac{1}{(2\pi i)^2} \oint_{\Sigma} \frac{ds}{s} \oint_{C_a} \frac{du}{u} 
    \MeijerG{1,1}{1,r+1}{0}{0, -\nu_1, \ldots, -\nu_r}{-\frac{x}{s}}  \\
    \times 
    \MeijerG{r,0}{0,r}{-}{\nu_1, \ldots, \nu_r}{\frac{y}{u}}  
    \left( \prod_{j=1}^n \frac{s-a_j}{u-a_j} \right) \frac{1}{s-u},
  \end{multline}
  where $\Sigma$ is a closed, positive oriented contour around $0$, and $C_a$ is a closed contour, disjoint from $\Sigma$, in the right half-plane encircling each $a_j$ once in the positive direction.
\end{cor}

The above derivation of \eqref{eq:degenerate_Ginibre} was done under the assumption
that the $a_j$'s are mutually distinct. However in \eqref{eq:degenerate_Ginibre} we
can easily let some or all the $a_j$'s come together, and the expression is valid for every $A$
with non-zero squared singular values $a_1, \ldots, a_n$.

  In the limiting case where all $a_j \to 1$ and so $A = I$ we obtain
  \begin{multline} \label{eq:alternative_product_Ginibre}
    K_n(x,y) = 
    \frac{1}{(2\pi i)^2} \oint_{\Sigma} \frac{ds}{s} \oint_{\Gamma} \frac{du}{u} 
    \MeijerG{1,1}{1,r+1}{0}{0, -\nu_1, \ldots, -\nu_r}{-\frac{x}{s}} \\
    \times
    \MeijerG{r,0}{0,r}{-}{\nu_1, \ldots, \nu_r}{\frac{y}{u}} 
    \frac{(s-1)^n}{(u-1)^n}  \frac{1}{s-u}
  \end{multline}
  with $\Sigma$ and $\Gamma$ disjoint closed contours where $\Sigma$ encloses $0$ and $\Gamma$
	encloses~$1$. 
	This is an alternative expression for \eqref{kernel product}, which 
	for $r=1$ reduces to \eqref{eq:traditional_Wishart}.

\section{Multiplication with truncated unitary matrices}\label{section: trunc}


\subsection{Proof of Theorem \ref{TruncationTransformationKernel}}\label{section: proof truncated}

Theorem \ref{TruncationTransformationKernel} follows from Lemma \ref{Transformation-phi} 
in the same way as Theorem \ref{GinibreTransformationKernel} does.

In the situation of Theorem \ref{TruncationTransformationKernel} we have
\begin{equation}
  \varphi(t) =
  \begin{cases}
    t^{\nu}(1-t)^{\mu-1} & \text{ for } 0 < t < 1, \\
    0 & \text{ otherwise,}
  \end{cases}
\end{equation}
with integers $\nu \geq 0$ and $\mu \geq 1$. Then
\begin{equation}
  b_j^{-1} = \int_0^1 t^{j+\nu} (1-t)^{\mu-1} dt =   \frac{(\mu-1)!(j+\nu)!}{(j+\nu + \mu)!}.
\end{equation}

We take for $\psi$
\begin{equation}
  \begin{split}
    \psi(x) = {}& \sum_{j=-\nu}^{\infty} b_j x^j
    = \frac{1}{(\mu-1)!} \sum_{j=-\nu}^{\infty} \frac{(j+\nu+\mu)!}{(j+\nu)!} x^j \\
    = {}& \mu x^{-\nu} (1-x)^{-\mu-1},
    \qquad \text{ for } 0 < |x| < 1. 
  \end{split}
\end{equation}

Then the statements of Theorem \ref{TruncationTransformationKernel} follow immediately
from Lemma \ref{Transformation-phi}. \qed

\subsection{Proof of Corollary \ref{cor:truncated_iterated}}\label{section: proof corollary trunc}

The proof is similar to that of Corollary \ref{cor:Ginibre_iterated}. Similar to the function $\varphi$ in \eqref{eq:varphi_multiple_G}, we take here
\begin{equation}
  \varphi = \varphi_r \ast \cdots \ast \varphi_1
\end{equation}
with now
\begin{equation}
  \varphi_k(t) =
  \begin{cases}
    \frac{1}{
    (\mu_k-1)!} t^{\nu_k} (1-t)^{\mu_k-1}, & \text{ for } 0 < t < 1 \\
    0 & \text{ otherwise}.
  \end{cases}
\end{equation}
Then $\varphi$ is supported on $[0,1]$ with moments
\begin{equation}
  b_j^{-1} = \int_0^1 t^j \varphi (t) dt = 
  \prod_{k=1}^r \frac{(j+\nu_k)!}{(j+\nu_k+\mu_k)!}.
\end{equation}

We can write $\varphi$ as a Meijer G-function
\begin{equation} \label{Truncationphi} 
	\varphi(t) = \MeijerG{r,0}{r,r}{\nu_1+\mu_1, \ldots, \nu_r + \mu_r}{\nu_1, \ldots, \nu_r}{t}. 
	\end{equation}
Furthermore, analogous to \eqref{psi:defseries}, we take here
\begin{equation}
  \psi(s) = \sum_{j=0}^{\infty} b_j s^j = \sum_{j=0}^{\infty} \left( \prod_{k=1}^r \frac{(j+\nu_k+ \mu_k)!}{(j+\nu_k)!} \right) s^j,
\end{equation}
which is a generalized hypergeometric series
\begin{equation}
  \psi(s) = \prod_{k=1}^r	\frac{(\nu_k+\mu_k)!}{\nu_k!}
  \hyperF{r + 1}{r}{1, \nu_1 + \mu_1 + 1, \dotsc, \nu_r + \mu_r + 1}{\nu_1 + 1, \dotsc, \nu_r + 1}{s}
\end{equation}
and also a Meijer G-function
\begin{equation} \label{Truncationpsi}
  \psi(s) = \MeijerG{1,r+1}{r+1,r+1}{0,-\nu_1-\mu_1, \ldots, -\nu_r-\mu_r}{0,-\nu_1, \ldots, -\nu_r}{-s}.
\end{equation}
Then Corollary \ref{cor:truncated_iterated} is a direct consequence of Lemma \ref{Transformation-phi}\ref{enu:ransformation-phi_c}.
\qed

\subsection{Proof of Theorem \ref{TruncationProduct}}
\label{section: proof TruncationProduct}

\subsubsection{Proof of the necessity of \eqref{TruncationProductCondition} }

The necessity of the condition \eqref{TruncationProductCondition} follows
from the following linear algebra lemma.
\begin{lemma} \label{lem:linear_algebraic}
Suppose that 
\begin{equation} \label{OppositeCondition} 
	d = n - \sum_{j=1}^r \mu_j \geq 1. 
	\end{equation}
Then the product $Y_r = T_r \cdots T_1$, where $T_1, \dotsc, T_r$ are defined in Assumption \ref{defn:T_matrices} with $\nu_0 = 0$, has a singular value at $1$ of multiplicity $\geq d$.
\end{lemma}
Note that although $T_1, \dotsc, T_r$ are random, Lemma \ref{lem:linear_algebraic} uses only their matrix structure and the result and argument are deterministic.

\begin{proof}
Each $T_j$ is the truncation of a unitary matrix $U_j$, which we partition as
\begin{equation}
  U_j =
  \begin{pmatrix}
    T_j & * \\
    S_j & *
  \end{pmatrix}
\end{equation}
where $S_j$ has size $\mu_j \times (n+\nu_{j-1})$, and $*$ denotes a block that is not important for our present purpose. Thus
\begin{equation} \label{rankSk} 
	\rank(S_j) \leq \mu_j. 
	\end{equation}
Note that 
\begin{equation}
  U_j^* U_j = \begin{pmatrix}
    T_j^* T_j + S_j^* S_j & * \\
    * & *
  \end{pmatrix}
\end{equation}
and therefore, since $U_j$ is unitary,
\begin{equation} \label{eq:TkSk} 
	T_j^* T_j +  S_j^* S_j = I.
	\end{equation}

Below we denote $Y_j = T_j T_{j - 1} \dotsm T_1$. Taking $j = 1$, we find by \eqref{rankSk} and \eqref{eq:TkSk} that
\begin{equation} \label{rankY1}
	\rank(I-Y_1^* Y_1) = \rank(S_1^*S_1) \leq \mu_1.
	\end{equation}
For $1 \leq j < r$ we have $Y_{j+1} = T_{j + 1} Y_j$ and therefore by \eqref{eq:TkSk}
\begin{equation}
  I - Y_{j+1}^* Y_{j+1} = I - Y_j^* T_{j + 1}^* T_{j + 1} Y_j  = I - Y_j^* Y_j + Y_j^* S_{j + 1}^* S_{j + 1} Y_j.
\end{equation}
Using \eqref{rankSk} and elementary properties of the rank, we then get from this	
\begin{equation} \label{rankYk} 
  \rank(I-Y_{j+1}^* Y_{j+1}) \leq \rank(I-Y_j^* Y_j) + \mu_{j + 1}.
\end{equation}
The two inequalities \eqref{rankY1} and \eqref{rankYk}, and the assumption
\eqref{OppositeCondition} of the lemma, lead to
\begin{equation} \label{rankYr}
  \rank(I - Y_r^* Y_r) \leq \sum_{j=1}^r \mu_j  = n-d
\end{equation}
with $d \geq 1$. Since $Y_r^*Y_r$ has size $n \times n$, we obtain
\begin{equation}
  \dim \kernel(I-Y_r^* Y_r) \geq d.
\end{equation}
Every vector $v$ in the kernel of $I-Y_r^* Y_r$ is a right singular vector of $Y_r$
with singular value $1$. Thus $Y_r$ has at least $d$ singular values at $1$.
\end{proof}

\subsubsection{Proof of \eqref{Pk-ProductTruncation}, \eqref{Qk-ProductTruncation}, and 
\eqref{Kn-ProductTruncation} under condition \eqref{TruncationProductCondition}}

First we consider a limiting degenerate case of Corollary \ref{cor:truncated_iterated}, 
with $\nu_0 = 0$ and the matrix $X$ replaced by the deterministic matrix $A$ defined in \eqref{eq:defn_A}. 
Then similar to Corollary \ref{cor:degenerate_Ginibre} that is proved in 
Section \ref{subsec:proof_degenerate_Ginibre}, the biorthogonal system and the correlation kernel 
for the squared singular values of $T_r \dotsm T_1 A$ are expressed by \eqref{TruncationPk}, 
\eqref{TruncationQk} and \eqref{trancationK_n} with the degenerate form of $q_k$ and $k_n$
containing Dirac $\delta$-functions. 
The $p_k$'s are polynomials and they are not degenerate.

Here we take the biorthogonal system \eqref{eq:pkdegenerate} and \eqref{eq:qkdegenerate} for 
the degenerate polynomial ensemble of the squared singular values of $A$. Then since 
$\int_0^{\infty} \varphi(t) \delta(y/t-a_k) dt/t = \varphi(y/a_k)/a_k$, it follows from \eqref{TruncationQk}, that 
\begin{equation} \label{Qkwithaj} 
	Q_k(y) 
 = \frac{1}{\Delta_{k+1}(a)} 
	\begin{vmatrix} 1 & 1 & \cdots & 1 \\
		a_1 & a_2 &  \cdots & a_{k+1} \\
		\vdots & \vdots & \ddots & \vdots \\
		\varphi(y/a_1)/a_1 & \varphi(y/a_2)/a_2  & \cdots & \varphi(y/a_{k+1})/a_{k+1} \end{vmatrix}
		\end{equation}
where $\varphi$ is defined in \eqref{Truncationphi}. On the other hand, $P_k$ 
is obtained simply by plugging \eqref{eq:pkdegenerate} into \eqref{TruncationPk}.

Now we take the limit where $A \to I$, so that all $a_j$ tend to $1$.
Then $P_k$ becomes the Hadamard product of \eqref{Truncationpsi}
with 
\begin{equation}
  p_k(x) = (x-1)^k = \sum_{j=0}^k (-1)^{k-j} \binom{k}{j} x^j
\end{equation}
which is 
\begin{equation}
  P_k(x) = \sum_{j=0}^k (-1)^{k-j} \binom{k}{j} b_j x^j, \qquad 
  b_j = \prod_{l=1}^r \frac{(j+\nu_l+\mu_l)!}{(j+\nu_l)!}
  \end{equation}
as also given in \eqref{Pk-ProductTruncation}.

The limit $a_j \to 1$ in \eqref{Qkwithaj} has to be done with more care.
Note that each $\varphi_j(x)$ is zero for $x > 1$, real analytic on $(0,1)$,
with the behaviour $\varphi_j(x) \sim \Gamma(\mu_j)^{-1} (1 - x)^{\mu_j - 1}$ as $x \to 1_-$. Then inductively, for $\varphi^{(j)} = \varphi_j * \varphi_{j - 1} * \dotsm * \varphi_1$, it is zero for $x > 1$ and 
\begin{equation}
  \begin{aligned}[b]
    \varphi^{(j)}(x) = {}& \varphi_j * (\varphi_{j - 1} * \dotsm * \varphi_1)(x) \\
    = {}& \int^{\infty}_0 \varphi_j \left( \frac{x}{t} \right) (\varphi_{j - 1} * \dotsm * \varphi_1)(t) \frac{dt}{t} \\
    \sim {}& c_j\int^1_x \left( 1 - \frac{x}{t} \right)^{\mu_j - 1} (1 - t)^{\mu_1 + \cdots + \mu_{j - 1} - 1} \frac{dt}{t} \\
    \sim {}& c_j'(1 - x)^{\mu_1 + \cdots + \mu_j - 1} ,
  \end{aligned}
   \quad \text{ as } x \to 1_-,
\end{equation}
for some non-zero constants $c_j, c_j'$.
Thus $\varphi = \varphi^{(r)}$ has only $\mu_1 + \cdots + \mu_r-1$ derivatives on $(0,\infty)$.

However, since $n \leq \sum_{j=1}^r \mu_j$ and $k \leq n-1$,
we have that $\varphi$ is $k$ times differentiable on $(0,\infty)$ and
then it follows that the limit of \eqref{TruncationQk} for all $a_j \to 1$
can be evaluated as 
\begin{equation} 
  \begin{split}
    Q_k(y) = {}& \frac{1}{k!}  \left. \left(\frac{d}{da}\right)^k \right|_{a=1}  
		\left[\frac{1}{a} \varphi\left( \frac{y}{a}\right) \right]  \\
    = {}& \frac{(-1)^k}{k!} \left(\frac{d}{dy} \right)^k \left[ y^k \varphi(y)) \right].
  \end{split}
\end{equation}
The second equality can be proved by induction on $k$.
Recall that $\varphi$ is the Meijer G-function \eqref{Truncationphi}
and then by elementary properties of Meijer G-functions, we find
that $Q_k$ is given by \eqref{Qk-ProductTruncation}.

Then \eqref{Pk-ProductTruncation} and \eqref{Qk-ProductTruncation} give
a biorthogonal system for the squared singular values of $T_r \cdots T_1$,
which follow therefore a polynomial ensemble.
Having the biorthogonal system \eqref{Pk-ProductTruncation} and \eqref{Qk-ProductTruncation} we can proceed as in the proof of 
Proposition 2.7 in \cite{Kieburg-Kuijlaars-Stivigny15} to find the
expression \eqref{Kn-ProductTruncation} for the correlation kernel.
\qed

\begin{rmk} \label{eq:truncated_KKS}

In case $\mu_1 \geq n$, we can also use Theorem \ref{TruncationTransformationKernel}
to prove  \eqref{Kn-ProductTruncation} by induction on $r$. The argument is similar to 
the proof of \eqref{kernel product} given in Section \ref{subsec:product_Ginibre_KZ}
and we do not give the details. We only mention that instead of  \eqref{GammaIntegral} and 
\eqref{GammaRecipIntegral} we now use 
\begin{align}
  \frac{\mu_r}{2\pi i} \int_L s^{-\nu_r-u-1} (1 - s)^{-\mu_r - 1} ds 
		& = \frac{\Gamma(u + \nu_r + \mu_r + 1)}{\Gamma(\mu_r) \Gamma(u + \nu_r + 1)}, \label{eq:NIST_needed} \\
	\int^1_0 t^{\nu_r+v} (1 - t)^{\mu_r - 1} dt & = \frac{\Gamma(\mu_r) \Gamma(v + \nu_r + 1)}{\Gamma(v + \nu_r + \mu_r + 1)},
	\label{BetaIntegral}
\end{align}
to obtain the desired result. Note that the integral in \eqref{eq:NIST_needed} 
is equivalent to \cite[5.12.9]{Boisvert-Clark-Lozier-Olver10}, while \eqref{BetaIntegral} is the familiar Beta integral.
\end{rmk}

\begin{rmk}
For $r=1$ the kernel is the correlation
kernel of a Jacobi unitary ensemble \cite[Section 3.8.3]{Forrester10}, \cite{Zyczkowski-Sommers00}. 
The double contour integral formula \eqref{Kn-ProductTruncation} with $r = 1$  for this kernel
is not well-known, and it was first described, to the best knowledge of the authors, 
in \cite[Theorem 3(d)]{Adler-van_Moerbeke-Wang11} as a special case of 
the \emph{Jacobi--\Pineiro\ minor process}.  (See also \cite{Forrester-Wang15} for its relation to the Jacobi Muttalib--Borodin model.) 

In the notation of \cite{Adler-van_Moerbeke-Wang11}, 
we take fixed time $n$, $\alpha_k = n + \nu_1 - k,$ for $k = 1, \dotsc, n$ and $M' = \mu_1$. 
Then
\begin{equation}
  \left. K_n(x, y) \right\rvert_{r = 1} = \left. \frac{y^{\nu_1}}{x^{\nu_1}} K(n, y; n, x) \right\rvert_{\text{$\alpha_k = n + \nu_1 - k$ for $k = 1, \dotsc, n$, $M' = \mu_1$}},
\end{equation}
analogous to \eqref{eq:K_n_r=1_equivalence}, where $K$ is the correlation kernel defined in \cite[Theorem 3(d)]{Adler-van_Moerbeke-Wang11}.
\end{rmk}


\end{document}